\documentclass[11pt,a4paper]{article}

\usepackage{theorem,enumerate}
\usepackage{amsmath,latexsym,amssymb,amsfonts}
\usepackage{eucal}
\usepackage{color}
\usepackage{comment}
\usepackage{cite}
\usepackage{graphics}

\usepackage{here}
\usepackage{fancybox}
\usepackage[all]{xy}
\usepackage{tabularx}
   \newcolumntype{C}{>{\centering\arraybackslash}X}
   \newcolumntype{L}{>{\raggedright\arraybackslash}X}
   \newcolumntype{R}{>{\raggedleft\arraybackslash}X}
   \newcolumntype{D}{>{\scriptsize\centering\arraybackslash}X}

\theorembodyfont{\normalfont\slshape}
\newtheorem{thm}{Theorem}[section]

\newtheorem{prop}[thm]{Proposition}
\newtheorem{lem}[thm]{Lemma}

\newtheorem{claim}{Claim}[section]

\newtheorem{conj}{Conjecture}

\newtheorem{remark}{Remark}

\numberwithin{equation}{section}

\newcommand{\proof}{\medbreak\noindent\textit{Proof.}\quad}

\newcommand{\qed}{{$\quad\square$\vs{3.6}}}

\newcommand{\vs}[1]{\vspace*{#1 mm}}

\def\AA{{ \mathcal{A}}}

\def\HH{{ \mathcal{H}}}

\def\PP{{ \mathcal{P}}}

\def\RR{{ \mathcal{R}}}
\def\SS{{ \mathcal{S}}}
\def\TT{{ \mathcal{T}}}

\newcommand{\fb}{\fbox}
\newcommand{\db}{\doublebox}

\bfseries\normalfont

\title{Ramsey-type problems on induced covers and induced partitions toward the Gy\'{a}rf\'{a}s-Sumner conjecture}

\author{
Shuya Chiba$^{1}$\footnote{\texttt{e-mail:schiba@kumamoto-u.ac.jp}}\and \
Michitaka Furuya$^{2}$\footnote{\texttt{e-mail:michitaka.furuya@gmail.com}}\vs{5}\\
$^{1}$\textsl{Applied Mathematics, Faculty of Advanced Science and Technology,}\\
\textsl{Kumamoto University,}\\
\textsl{2-39-1 Kurokami, Kumamoto 860-8555, Japan}\\
$^{2}$\textsl{College of Liberal Arts and Sciences,}\\
\textsl{Kitasato University,}\\
\textsl{1-15-1 Kitasato, Minami-ku, Sagamihara, Kanagawa 252-0373, Japan}
}

\date{}

\begin{document}

\maketitle

\begin{abstract}
Gy\'{a}rf\'{a}s and Sumner independently conjectured that for every tree $T$, there exists a function $f_{T}:\mathbb{N}\rightarrow \mathbb{N}$ such that every $T$-free graph $G$ satisfies $\chi (G)\leq f_{T}(\omega (G))$, where $\chi (G)$ and $\omega (G)$ are the {\it chromatic number} and the {\it clique number} of $G$, respectively.
This conjecture gives a solution of a Ramsey-type problem on the chromatic number.

For a graph $G$, the {\it induced SP-cover number ${\rm inspc}(G)$} (resp. the {\it induced SP-partition number ${\rm inspp}(G)$}) of $G$ is the minimum cardinality of a family $\mathcal{P}$ of induced subgraphs of $G$ such that each element of $\mathcal{P}$ is a star or a path and $\bigcup _{P\in \mathcal{P}}V(P)=V(G)$ (resp. $\dot\bigcup _{P\in \mathcal{P}}V(P)=V(G)$).
Such two invariants are directly related concepts to the chromatic number.
From the viewpoint of this fact, we focus on Ramsey-type problems for two invariants ${\rm inspc}$ and ${\rm inspp}$, which are analogies of the Gy\'{a}rf\'{a}s-Sumner conjecture, and settle them.
As a corollary of our results, we also settle other Ramsey-type problems for widely studied invariants.
\end{abstract}



\section{Introduction}\label{sec1}

All graphs considered in this paper are finite, simple and undirected.
For terms and symbols not defined in this paper, we refer the reader to \cite{D}.
For two graphs $G$ and $H$, $G$ is said to be {\it $H$-free} if $G$ contains no copy of $H$ as an induced subgraph.
For a family $\HH$ of graphs, a graph $G$ is said to be {\it $\HH$-free} if $G$ is $H$-free for every $H\in \HH$.
In this context, the members of $\HH$ are called {\it forbidden subgraphs}.
For two families $\HH_{1}$ and $\HH_{2}$ of graphs, we write $\HH_{1}\leq \HH_{2}$ if for every $H_{2}\in \HH_{2}$, $H_{2}$ contains a copy of an element of $\HH_{1}$ as an induced subgraph.
The relation ``$\leq $'' between two families of forbidden subgraphs was introduced in \cite{FKLOPS}.
Note that if $\HH_{1}\leq \HH_{2}$, then every $\HH_{1}$-free graph is also $\HH_{2}$-free.
Let $K_{n}$, $K_{1,n}$ and $P_{n}$ denote the {\it complete graph} of order $n$, the {\it star} of order $n+1$ and the {\it path} of order $n$, respectively.
Let $\overline{K_{n}}$ denote the complement of $K_{n}$, i.e., the graph of order $n$ having no edges.

Gy\'{a}rf\'{a}s~\cite{G} and Sumner~\cite{S} independently conjectured that for every tree $T$, there exists a function $f_{T}:\mathbb{N}\rightarrow \mathbb{N}$ such that every $T$-free graph $G$ satisfies $\chi (G)\leq f_{T}(\omega (G))$, where $\chi (G)$ and $\omega (G)$ are the {\it chromatic number} and the {\it clique number} of $G$, respectively.
In graph theory, the Gy\'{a}rf\'{a}s-Sumner conjecture has been studied as a fundamental and important problem.
The conjecture is true for the case where $T$ is a path, a star, and some specific trees, but it is widely open (see a survey~\cite{SS}).

Here we focus on Ramsey-type problems.
A classical Ramsey theorem~\cite{R} is equivalent to the following statement:
For a family $\HH$ of graphs, there exists a constant $c(\HH)$ such that $|V(G)|\leq c(\HH)$ for every $\HH$-free graph $G$ if and only if $\HH\leq \{K_{n},\overline{K_{n}}\}$ for an integer $n\geq 1$.
Furthermore, an analogy of the Ramsey theorem for connected graphs is also well-known as a folklore (see \cite[Proposition~9.4.1]{D}):
For a family $\HH$ of connected graphs, there exists a constant $c(\HH)$ such that $|V(G)|\leq c(\HH)$ for every connected $\HH$-free graph $G$ if and only if $\HH\leq \{K_{n},K_{1,n},P_{n}\}$ for an integer $n\geq 1$.
Remark that the $\overline{K_{n}}$-freeness and the $\{K_{1,n},P_{n}\}$-freeness play a similar role in above results.
In particular, the combination of induced stars and induced paths is a natural alternative concept of an independent set when we treat induced subgraphs.

Recently, similar type problems on invariants (other than the order) were obtained as Ramsey-type problems (see, for example, \cite{BPRS,CF,CFKP,DDL,F}).
For a graph-invariant $\mu $, we consider the following condition concerning a (finite) family $\HH$ of graphs:
\begin{enumerate}[{\bf (P-$\mu $)}]
\item
There exists a constant $c(\HH)$ such that $\mu (G)\leq c(\HH)$ for every connected $\HH$-free graph $G$.
\end{enumerate}
A Ramsey-type problem of $\mu $ is to characterize $\HH$ satisfying {\rm (P-$\mu $)}.

Since there exist infinitely many graphs with sufficiently large chromatic number and sufficiently large girth~\cite{E}, we can easily verify that the following conjecture is equivalent to the Gy\'{a}rf\'{a}s-Sumner conjecture.

\begin{conj}
\label{GSconj}
Let $\HH$ be a finite family of connected graphs.
Then $\HH$ satisfies {\rm (P-$\chi $)} if and only if $\HH\leq \{K_{n},T\}$ for an integer $n\geq 1$ and a tree $T$.
\end{conj}

Thus, settling Ramsey-type problems for invariants closely related to the chromatic number, one might expect to obtain some essential information and effective techniques for the Gy\'{a}rf\'{a}s-Sumner conjecture.
As a study of this policy, for example, Chudnovsky and Seymour~\cite{CS} posed a conjecture on a Ramsey-type problem for the cochromatic number, where the {\it cochromatic number} of a graph $G$ is the minimum cardinality of a family of independent sets and cliques of $G$ whose union is $V(G)$, and proved that the conjecture is equivalent to the Gy\'{a}rf\'{a}s-Sumner conjecture.
On the other hand, as far as we know, no Ramsey-type problem for an invariant around the chromatic number was completely settled.
Our main aim in this paper is to approach the Gy\'{a}rf\'{a}s-Sumner conjecture by settling such Ramsey-type problems.
Remark that Choi, Kim and Park~\cite{CKP} characterized forbidden structures which are not necessarily induced, bounding some analogies of the chromatic number by a constant.

Let $G$ be a graph, and let $\AA$ be a family of graphs.
A family $\PP$ of induced subgraphs of $G$ is called an {\it induced $\AA$-cover} of $G$ if $\bigcup _{P\in \PP}V(P)=V(G)$ and each element of $\PP$ is isomorphic to an element of $\AA$.
Note that some elements of an induced $\AA$-cover of $G$ might have common vertices.
An induced $\AA$-cover $\PP$ of $G$ is called an {\it induced $\AA$-partition} of $G$ if the elements of $\PP$ are pairwise vertex-disjoint.
We mainly consider the case where $\AA$ is equal to the family $\AA_{\rm sp}:=\{K_{1,n}:n\geq 1\}\cup \{P_{n}:n\geq 1\}$, i.e., the family of stars and paths.
An induced $\AA_{\rm sp}$-cover (resp. an induced $\AA_{\rm sp}$-partition) of $G$ is called an {\it induced SP-cover} (resp. an {\it induced SP-partition}) of $G$.
The minimum cardinality of an induced SP-cover (resp. an induced SP-partition) of $G$, denoted by ${\rm inspc}(G)$ (resp. ${\rm inspp}(G)$), is called the {\it induced SP-cover number} (resp. the {\it induced SP-partition number}) of $G$.

Let us now compare the induced SP-cover/partition numbers and the chromatic number.
Note that the chromatic number of a graph $G$ is equal to the minimum cardinality of an induced $\{\overline{K_{n}}:n\geq 1\}$-cover of $G$, which is also equal to the minimum cardinality of an induced $\{\overline{K_{n}}:n\geq 1\}$-partition of $G$.
As we mentioned in the third paragraph of this section, the combination of induced stars and induced paths is related to the existence of independent sets, or induced $\overline{K_{n}}$'s.
In fact, for an induced SP-cover $\PP$ of $G$, since every element of $\PP$ is a bipartite graph, $\{G[A_{P}],G[B_{P}]:P\in \PP\}$ is an induced $\{\overline{K_{n}}:n\geq 1\}$-cover of $G$, where $\{A_{P},B_{P}\}$ is a bipartition of $P$.
Hence we obtain $\chi (G)\leq 2 \, {\rm inspc}(G)\leq 2 \, {\rm inspp}(G)$ for every graph $G$, which implies that the induced SP-cover/partition numbers are directly related concepts to the chromatic number.
In particular, if we restrict induced $\AA$-covers/partitions to the case where each element of $\AA$ is connected, then small induced SP-covers/partitions are alternatives to small induced $\{\overline{K_{n}}:n\geq 1\}$-covers/partitions, and in this sense those two numbers can be regarded as ``connected versions'' of the chromatic number.
In this paper, we characterize the finite families $\HH$ of connected graphs satisfying {\rm (P-${\rm inspc}$)} or {\rm (P-${\rm inspp}$)}.

\subsection{Fundamental definitions}\label{sec1.1}

Let $G$ be a graph.
Let $V(G)$ and $E(G)$ denote the {\it vertex set} and the {\it edge set} of $G$, respectively.
For a vertex $x\in V(G)$, let $N_{G}(x)$ denote the {\it neighborhood} of $x$ in $G$, i.e., $N_{G}(x)=\{y\in V(G): xy\in E(G)\}$.
For two vertices $x,y\in V(G)$, let ${\rm dist}_{G}(x,y)$ denote the length of a shortest $x$-$y$ path of $G$.
For a vertex $x\in V(G)$, let ${\rm ecc}_{G}(x)=\max\{{\rm dist}_{G}(x,y):y\in V(G)\}$.
The {\it diameter} of $G$, denoted by ${\rm diam}(G)$, is defined as ${\rm diam}(G)=\max\{{\rm ecc}_{G}(x):x\in V(G)\}$.
Let $\alpha (G)$ denote the {\it independence number} of $G$, i.e., the maximum cardinality of an independent set of $G$.
For a subset $X$ of $V(G)$, let $G[X]$ (resp. $G-X$) denote the subgraph of $G$ induced by $X$ (resp. $V(G)\setminus X$).
For a subset $F$ of $E(G)$, let $G-F$ denote the spanning subgraph of $G$ with the edge set $E(G)\setminus F$.
For two integers $n_{1}\geq 1$ and $n_{2}\geq 1$, the {\it Ramsey number} $R(n_{1},n_{2})$ is the minimum positive integer $R$ such that any graph of order at least $R$ contains a clique of cardinality $n_{1}$ or an independent set of cardinality $n_{2}$.

Recall that we have defined the family $\AA_{\rm sp}$.
In this paper, we additionally focus on induced $\AA$-covers/partitions of graphs for the case where $\AA$ is equal to one of the following families:
\begin{enumerate}[{$\bullet $}]
\setlength{\parskip}{0cm}
\setlength{\itemsep}{0cm}
\item
$\AA_{\rm s}:=\{K_{1}\}\cup \{K_{1,n}:n\geq 1\}$, i.e., the family of stars where the graph of order one is regarded as a star.
\item
$\AA_{\rm p}:=\{P_{n}:n\geq 1\}$, i.e., the family of paths.
\end{enumerate}
Then for a graph $G$, we can naturally define some terminologies and notations as follows:
\begin{enumerate}[{$\bullet $}]
\setlength{\parskip}{0cm}
\setlength{\itemsep}{0cm}
\item
An induced $\AA_{\rm s}$-cover (resp. an induced $\AA_{\rm s}$-partition) of $G$ is called an {\it induced star cover} (resp. an {\it induced star partition}) of $G$.
\item
The minimum cardinality of an induced star cover (resp. an induced star partition) of $G$, denoted by ${\rm insc}(G)$ (resp. ${\rm insp}(G)$), is called the {\it induced star cover number} (resp. the {\it induced star partition number}) of $G$.
\item
An induced $\AA_{\rm p}$-cover (resp. an induced $\AA_{\rm p}$-partition) of $G$ is called an {\it induced path cover} (resp. an {\it induced path partition}) of $G$.
\item
The minimum cardinality of an induced path cover (resp. an induced path partition) of $G$, denoted by ${\rm inpc}(G)$ (resp. ${\rm inpp}(G)$), is called the {\it induced path cover number} (resp. the {\it induced path partition number}) of $G$.
\end{enumerate}
The induced path cover/partition problems themselves are interesting topics in graph theory, and have been studied (see \cite{AHW,CMSHH,HSMW,HSMW2,PC3}).
Furthermore, research on the induced star cover/partition problems have recently been initiated (see \cite{DV,MV,SVSM}).
The induced path partition number is frequently called the {\it induced path number}.
However, in this paper, we use the terminology ``induced path partition number'' to avoid confusion.

As we mentioned above, our main aim is to characterize the finite families $\HH$ of connected graphs satisfying {\rm (P-${\rm inspc}$)} or {\rm (P-${\rm inspp}$)}.
Furthermore, during the process, we also discuss about the conditions {\rm (P-${\rm insc}$)}, {\rm (P-${\rm insp}$)}, {\rm (P-${\rm inpc}$)} and {\rm (P-${\rm inpp}$)}, and characterize the finite families satisfying one of them.
Remark that the authors~\cite{CF,CF2} recently settled Ramsey-type problems for the (not necessarily induced) path cover/partition numbers of graphs and digraphs.

\subsection{Main results}\label{sec1.2}

\begin{figure}
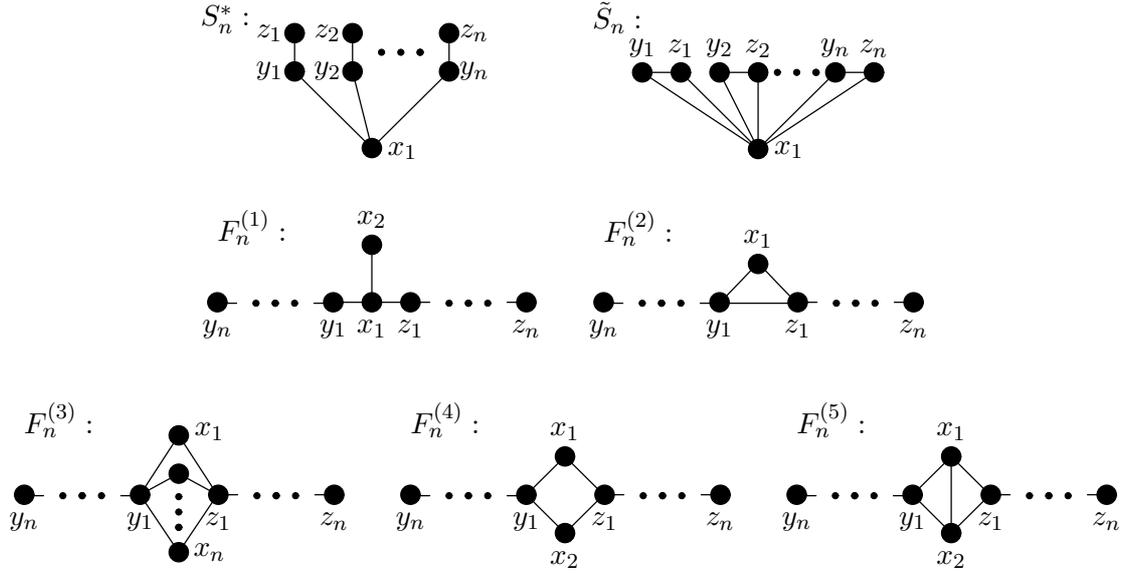

\begin{center}
{\unitlength 0.1in%
}%

\caption{Graphs $S^{*}_{n}$, $\tilde{S}_{n}$ and $F^{(i)}_{n}~(1\leq i\leq 5)$}
\label{f1}
\end{center}
\end{figure}

To state our main results, we prepare some graphs which will be used as forbidden subgraphs (see Figure~\ref{f1}).
Let $n\geq 2$ be an integer, and let $A:=\{x_{i},y_{i},z_{i}:1\leq i\leq n\}$.
\begin{enumerate}[{$\bullet $}]
\setlength{\parskip}{0cm}
\setlength{\itemsep}{0cm}
\item
Let $S^{*}_{n}$ denote the graph on $A\setminus \{x_{i}:2\leq i\leq n\}$ such that $E(S^{*}_{n})=\{x_{1}y_{i},y_{i}z_{i}:1\leq i\leq n\}$.
\item
Let $\tilde{S}_{n}$ be the graph obtained from $S^{*}_{n}$ by adding $n$ edges $x_{1}z_{i}~(1\leq i\leq n)$.
\item
Let $F^{(1)}_{n}$ denote the graph on $A\setminus \{x_{i}:3\leq i\leq n\}$ such that $E(F^{(1)}_{n})=\{x_{1}x_{2},x_{1}y_{1},x_{1}z_{1}\}\cup \{y_{i}y_{i+1},z_{i}z_{i+1}:1\leq i\leq n-1\}$.
\item
Let $F^{(2)}_{n}$ denote the graph on $A\setminus \{x_{i}:2\leq i\leq n\}$ such that $E(F^{(2)}_{n})=\{x_{1}y_{1},x_{1}z_{1},y_{1}z_{1}\}\cup \{y_{i}y_{i+1},z_{i}z_{i+1}:1\leq i\leq n-1\}$.
\item
Let $F^{(3)}_{n}$ denote the graph on $A$ such that $E(F^{(3)}_{n})=\{x_{i}y_{1},x_{i}z_{1}:1\leq i\leq n\}\cup \{y_{i}y_{i+1},z_{i}z_{i+1}:1\leq i\leq n-1\}$.
\item
Let $F^{(4)}_{n}$ be the graph obtained from $F^{(3)}_{n}$ by deleting $n-2$ vertices $x_{i}~(3\leq i\leq n)$.
\item
Let $F^{(5)}_{n}$ be the graph obtained from $F^{(4)}_{n}$ by adding the edge $x_{1}x_{2}$.
\end{enumerate}

Our main result is the following.

\begin{thm}
\label{mainthm}
Let $\HH$ be a finite family of connected graphs.
\begin{enumerate}[{\upshape(i)}]
\item
The family $\HH$ satisfies {\rm (P-${\rm inspc}$)} if and only if $\HH\leq \{K_{n},S^{*}_{n},F^{(1)}_{n},F^{(2)}_{n},F^{(3)}_{n}\}$ for an integer $n\geq 4$.
\item
The family $\HH$ satisfies {\rm (P-${\rm inspp}$)} if and only if $\HH\leq \{K_{n},S^{*}_{n},\tilde{S}_{n},F^{(1)}_{n},F^{(2)}_{n},F^{(4)}_{n},F^{(5)}_{n}\}$ for an integer $n\geq 4$.
\end{enumerate}
\end{thm}

\begin{figure}[H]
{
\footnotesize
\begin{center}
\[
\xymatrix@C=10pt@R=-1.4pt{
\ar@{--}[rrrrrr]\ar@{--}'[dd]&&&&&&\ar@{--}'[dd]\\
&\fb{Lem.~\ref{lem-2-indstar-01}}&&\fb{Lem.~\ref{lem-path-chibounded}}&&\db{Lem.~\ref{lem-2-indstar}}&\\
\ar@{--}[rrrrrr]&&&\ar@{=>}[dddd]&&&\\
&&&&&&\\
&&&&&&\\
&&&&&&\\
&&&\db{Prop.~\ref{prop-indstar}}\ar@{=>}[llddd]&&&\\
&&&&&&\\
\ar@{--}[rr]\ar@{--}'[dddd]&&\ar@{--}'[dddd]&&\ar@{--}[rr]\ar@{--}'[dddd]&&\ar@{--}'[dddd]\\
&\db{Prop.~\ref{mainthm-2}}\ar@{=>}[rrrr]&&&&\db{Prop.~\ref{mainthm-star/path}}&\\
&\db{Lem.~\ref{lem-nec-complete/star}}&&&&\db{Lem.~\ref{lem-nec-complete/star}}&\\
&\db{Lem.~\ref{lem-nec-pathtype}}&&&&\db{Lem.~\ref{lem-nec-pathtype}}&\\
\ar@{--}[rr]&\ar@{=>}[dddd]&&&\ar@{--}[rr]\ar@{=>}[ldddd]&\ar@{=>}[dddd]^(.4){\textup{\tiny{Rem.~\ref{remark-mainthm-2}}}}&\\
&&&&&&\\
&&&&&&\\
&&&&&&\\
&\underset{{\rm (inspc, inspp)}}{\db{Thm.~\ref{mainthm}}}&&\underset{{\rm (insc, insp, inpc, inpp)}}{\db{Thm.~\ref{mainthm-cor}}}&&{\db{Thm.~\ref{mainthm-3}}}&
}
\]
\end{center}
\vspace{-12pt}
\caption{The overall flowchart for results (the double box denotes our result)}
\label{overall flow}
}
\end{figure}
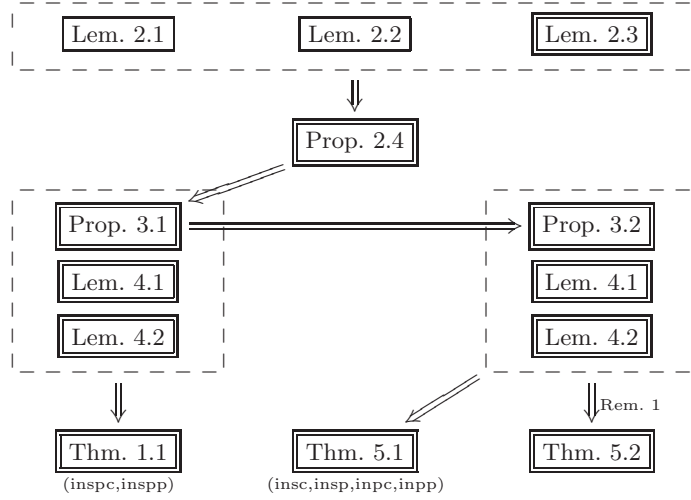

In Section~\ref{sec2}, we find small induced star covers/partitions of graphs with some forbidden subgraph conditions (Proposition~\ref{prop-indstar}).
In Section~\ref{sec3}, we essentially prove the ``if'' part of Theorem~\ref{mainthm} (Proposition~\ref{mainthm-2}), and as its corollary, we also obtain a result on induced star/path covers/partitions of graphs (Proposition~\ref{mainthm-star/path}).
In Section~\ref{sec-nec}, we calculate induced SP-cover/SP-partition numbers of some graphs whose constructions are important for the ``only if'' part of Theorem~\ref{mainthm} (Lemmas~\ref{lem-nec-complete/star} and \ref{lem-nec-pathtype}).
In Section~\ref{sec-proof}, we complete the proof of Theorem~\ref{mainthm}, and settle Ramsey-type problems for some invariants (Theorems~\ref{mainthm-cor} and \ref{mainthm-3}).
See also Figure~\ref{overall flow}.

\section{Induced star covers and induced star partitions}\label{sec2}

In this section, we give forbidden subgraph conditions which force the induced star cover/partition numbers to be small.

Let $G$ be a graph.
For subsets $X$ and $Y$ of $V(G)$, $X$ {\it dominates} $Y$ in $G$ if each vertex in $Y\setminus X$ is adjacent to a vertex in $X$.
A subset $X$ of $V(G)$ is a {\it dominating set} of $G$ if $X$ dominates $V(G)$ in $G$.
Let $K^{*}_{n}$ be the graph obtained from $K_{n}$ by adding a pendant edge to each vertex.
We recursively define the value $\alpha _{n,h}$ for integers $n\geq 1$ and $h\geq 1$ as follows:
$$
\begin{cases}
\alpha _{n,1}=1\\
\alpha _{n,h}=R(n,(n-1)\alpha _{n,h-1}+1)-1~~~(h\geq 2).
\end{cases}
$$
In \cite{F}, Furuya settled a Ramsey-type problem for the minimum cardinality of a dominating set of a graph, and he implicitly proved the following lemma.

\begin{lem}[Furuya~\cite{F}]
\label{lem-2-indstar-01}
Let $n\geq 1$ and $l\geq 1$ be integers.
Let $H$ be a connected $\{K^{*}_{n},S^{*}_{n}\}$-free graph with ${\rm diam}(H)\leq l$.
Then there exists a dominating set $U$ of $H$ with $|U|\leq R(n,n)\sum _{2\leq h\leq l}\alpha _{n,h}+1$.
\end{lem}


Our aim in this section is to give analogies of Lemma~\ref{lem-2-indstar-01} for the induced star cover/partition numbers (Proposition~\ref{prop-indstar}).

The following lemma is a special case of the results in \cite{CSS}.

\begin{lem}[Chudnovsky, Scott and Seymour~\cite{CSS}]
\label{lem-path-chibounded}
Let $n\geq 1$ be an integer.
Then there exists a constant $c_{\chi }(n)$ depending on $n$ only such that $\chi (G)\leq c_{\chi }(n)$ for every $\{K_{n},F^{(1)}_{n}\}$-free graph $G$.
\end{lem}

For integers $n\geq 3$ and $i\geq 1$, we let $\xi _{n,i}=\frac{(R(n-1,n)-1)^{i}-1}{R(n-1,n)-2}$.
Then the following lemma holds\footnote{The original proof of Lemma~\ref{lem-2-indstar} was long, but on the advice of Professor Jin Sun, based on \cite{Sun}, it could be made much shorter.}.

\begin{lem}
\label{lem-2-indstar}
Let $n$ and $m$ be integers with $n \geq m \geq 3$.
Let $G$ be a $\tilde{S}_{n}$-free graph, and let $x\in V(G)$ and $X\subseteq N_{G}(x)$. 
If $G[\{x\}\cup X]$ is $K_{m}$-free, then ${\rm insp}(G[\{x\}\cup X]) \leq \xi _{n, m-2}$.
\end{lem}
\proof
We prove it by induction on $m$. 
If $m = 3$, 
then 
$G[\{x\}\cup X]$ is an induced star; so ${\rm insp}(G[\{x\}\cup X]) = 1 = \xi _{n, 1}$. 
Thus we may assume that $m \ge 4$. 

Let $J$ be a maximal independent set of $G[X]$. 
We may assume that $X \setminus J \neq \emptyset$ 
(since otherwise, we have ${\rm insp}(G[\{x\}\cup X]) = 1 \le \xi _{n, m-2}$; 
so the assertion holds). 
Then $J$ dominates $X \setminus J$ in $G[X]$. 
We take a subset $I$ of $J$ dominating $X \setminus J$ in $G[X]$ 
so that $|I|$ is as small as possible. 
Set $I =\{y_{1}, \dots, y_{p}\}$. 
Since $I$ dominates $X \setminus J$ in $G[X]$, 
there is a partition $\{Y_{1}, \dots, Y_{p}\}$ of $X \setminus J$ 
such that $Y_{i} \subseteq N_{G}(y_{i})$ for $1 \le i \le p$. 
Then it follows from the minimality of $|I|$ that 
for each $i$ with $1 \le i \le p$, 
there is a vertex $z_{i}$ of $Y_{i}$ 
such that 
\begin{align}
\label{z_{i}y_{j}}
z_{i}y_{j} \notin E(G) \mbox{ for } 1 \le j \le p \mbox{ and } j \neq i. 
\end{align}
Let $I' = \{z_{i} : 1 \le i \le p\}$. 

Since $\{y_{i}\} \cup Y_{i} \subseteq X \subseteq N_{G}(x)$ 
for $1 \le i \le p$ 
and 
$G[\{x\} \cup X]$ is $K_{m}$-free, 
we see that $G[\{y_{i}\} \cup Y_{i}]$ is $K_{m-1}$-free 
for $1 \le i \le p$. 
Since $Y_{i} \subseteq N_{G}(y_{i})$, 
the induction hypothesis yields that 
\begin{align}
\label{y_{i} cup Y_{i}}
{\rm insp}(G[\{y_{i}\}\cup Y_{i}]) \leq \xi _{n, m-3} \mbox{ for } 1 \le i \le p. 
\end{align}

We show that
\begin{align}
\label{cond-lem-2-indstar-2}
p \leq R(n-1,n)-1.
\end{align}
Suppose that $p \geq R(n-1,n)$. 
Then 
$|I'| = |I| = p \geq R(n-1,n)$.
If $G[I']$ contains a copy of $K_{n-1}$, then 
$G[\{x\} \cup I']$ contains a copy of $K_{n}$; 
since $n \ge m$, $G[\{x\} \cup I']$ also contains a copy of $K_{m}$, 
which contradicts that $G[\{x\} \cup X]$ is $K_{m}$-free. 
Thus there exists a subset $\{z_{i_{1}}, z_{i_{2}}, \dots, z_{i_{n}}\}$ of $I'$ 
with $1 \le i_{1} < i_{2} < \cdots < i_{n} \le p$ 
such that $\{z_{i_{1}}, z_{i_{2}}, \dots, z_{i_{n}}\}$ is an independent set of $G$. 
Then by \eqref{z_{i}y_{j}}, 
we see that $G[\{x\} \cup \{y_{i_{\ell}}, z_{i_{\ell}} : 1 \le \ell \le n\}]$ 
is isomorphic to $\tilde{S}_{n}$, 
which contradicts that $G$ is $\tilde{S}_{n}$-free. 
Thus \eqref{cond-lem-2-indstar-2} is proved. 

Since $G[\{x\} \cup (J \setminus I)]$ is an induced star, 
it follows from \eqref{y_{i} cup Y_{i}} and \eqref{cond-lem-2-indstar-2} that 
\begin{align*}
{\rm insp}(G[\{x\}\cup X]) 
&\leq 
\left( \sum _{1 \le i \le p}{\rm insp}(G[\{y_{i}\} \cup Y_{i}]) \right) + 1\\
&\leq 
p \, \xi _{n, m - 3} + 1\\
&\le
\big( R(n-1,n) - 1 \big) \cdot \frac{\big( R(n-1,n) - 1 \big)^{m - 3} - 1}{R(n-1,n) - 2} \, + \, 1\\[2mm]
&= 
\frac{\big( R(n-1,n) - 1 \big)^{m - 2} - R(n-1,n) + 1}{R(n-1,n) - 2} \, + \, \frac{R(n-1,n) - 2}{R(n-1,n) - 2}\\[2mm]
&=
\frac{\big( R(n-1,n) - 1 \big)^{m - 2} - 1}{R(n-1,n) - 2} \\
&= 
\xi _{n, m - 2}. 
\end{align*}

This completes the proof of Lemma~\ref{lem-2-indstar}.
\qed

By combining Lemmas~\ref{lem-2-indstar-01}--\ref{lem-2-indstar}, we obtain the following proposition.

\begin{prop}
\label{prop-indstar}
Let $n\geq 3$ and $l_{0}\geq 1$ be integers.
Then there exist constants $c_{1}(n,l_{0})$ and $c_{2}(n,l_{0})$ depending on $n$ and $l_{0}$ only such that for a connected $\{K_{n},S^{*}_{n}\}$-free graph $H$ with ${\rm diam}(G)\leq l_{0}$,
\begin{enumerate}[{\upshape(1)}]
\item
if $H$ is $F^{(1)}_{n}$-free, then ${\rm insc}(H)\leq c_{1}(n,l_{0})$, and
\item
if $H$ is $\tilde{S}_{n}$-free, then ${\rm insp}(H)\leq c_{2}(n,l_{0})$.
\end{enumerate}
\end{prop}
\proof
We show that $c_{1}(n,l_{0})=(R(n,n)\sum _{2\leq h\leq l_{0}}\alpha _{n,h}+1)c_{\chi }(n)$ and $c_{2}(n,l_{0})=(R(n,n)\sum _{2\leq h\leq l_{0}}\alpha _{n,h}+1)\xi _{n,n-2}$ are desired constants, where $c_{\chi }(n)$ is the constant appearing in Lemma~\ref{lem-path-chibounded}.

Since $H$ is $K_{n}$-free, $H$ is also $K^{*}_{n}$-free.
Hence by Lemma~\ref{lem-2-indstar-01}, there exists a dominating set $U$ of $H$ with $|U|\leq R(n,n)\sum _{2\leq h\leq l_{0}}\alpha _{n,h}+1$.
For each vertex $v\in V(H)\setminus U$, take  a vertex $x_{v}\in N_{H}(v)\cap U$.
For a vertex $x\in U$, let $X_{x}=\{v\in V(H)\setminus U:x_{v}=x\}$.
Note that $X_{x}\subseteq N_{H}(x)$ for each $x\in U$, and $V(G)$ is the disjoint union of $\{x\}\cup X_{x}~(x\in U)$.
Hence for a fixed vertex $x\in U$, it suffices to prove the following two statements.
\begin{enumerate}
\setlength{\parskip}{0cm}
\setlength{\itemsep}{0cm}
\item[{\bf (S1)}]
If $H$ is $F^{(1)}_{n}$-free, then ${\rm insc}(H[\{x\}\cup X_{x}])\leq c_{\chi }(n)$.
\item[{\bf (S2)}]
If $H$ is $\tilde{S}_{n}$-free, then ${\rm insp}(H[\{x\}\cup X_{x}])\leq \xi _{n,n-2}$.
\end{enumerate}
By Lemma~\ref{lem-2-indstar}, (S2) clearly holds.

Now we prove that (S1) holds.
Suppose that $H$ is $F^{(1)}_{n}$-free.
Then by Lemmma~\ref{lem-path-chibounded}, $\chi (H[X_{x}])\leq c_{\chi }(n)$.
In particular, there exists an induced $\{\overline{K_{n}}:n\geq 1\}$-partition $\TT_{x}$ of $H[X_{x}]$ with $|\TT_{x}|\leq c_{\chi }(n)$.
Then $\SS_{x}:=\{H[\{x\}\cup T]:T\in \TT_{x}\}$ is an induced star cover of $H[\{x\}\cup X_{x}]$, and hence ${\rm insc}(H[\{x\}\cup X_{x}])\leq |\SS_{x}|=|\TT_{x}|\leq c_{\chi }(n)$, which leads to (S1).
\qed

\section{Induced SP-covers and induced SP-partitions}\label{sec3}

In this section, we prove the following proposition.

\begin{prop}
\label{mainthm-2}
Let $n\geq 4$ be an integer.
\begin{enumerate}[{\upshape(i)}]
\item
There exists a constant $c_{\rm inspc}(n)$ depending on $n$ only such that ${\rm inspc}(G)\leq c_{\rm inspc}(n)$ for every connected $\{K_{n},S^{*}_{n},F^{(1)}_{n},F^{(2)}_{n},F^{(3)}_{n}\}$-free graph $G$.
\item
There exists a constant $c_{\rm inspp}(n)$ depending on $n$ only such that ${\rm inspp}(G)\leq c_{\rm inspp}(n)$ for every connected $\{K_{n},S^{*}_{n},\tilde{S}_{n},F^{(1)}_{n},F^{(2)}_{n},F^{(4)}_{n},F^{(5)}_{n}\}$-free graph $G$.
\end{enumerate}
\end{prop}

\begin{remark}
\label{remark-mainthm-2}
A path $P$ of a graph $G$ is {\it isometric} if $P$ is a shortest path of $G$ joining the endvertices of $P$, i.e., ${\rm dist}_{G}(x,y)={\rm dist}_{P}(x,y)$ where $x$ and $y$ are the endvertices of $P$.
In the proof of Proposition~\ref{mainthm-2}, we construct small induced SP-covers/SP-partitions of a graph $G$ with some forbidden subgraph conditions.
One might notice that all paths belonging to such induced SP-covers/SP-partitions are also isometric paths of $G$.
By using the fact, we obtain solutions of Ramsey-type problems for some invariants concerning isometric paths (see Theorem~\ref{mainthm-3}).
\end{remark}

\medbreak\noindent\textit{Proof of Proposition~\ref{mainthm-2}.}\quad
Let $G$ be a connected $\{K_{n},S^{*}_{n},F^{(1)}_{n},F^{(2)}_{n},F^{(3)}_{n}\}$-free graph.
Since $F^{(4)}_{n}$ is an induced subgraph of $F^{(3)}_{n}$, it suffices to prove the following two statements.
\begin{enumerate}
\setlength{\parskip}{0cm}
\setlength{\itemsep}{0cm}
\item[{\bf (B1)}]
The value ${\rm inspc}(G)$ is bounded by a constant depending on $n$ only.
\item[{\bf (B2)}]
If $G$ is $\{\tilde{S}_{n},F^{(4)}_{n},F^{(5)}_{n}\}$-free, then ${\rm inspp}(G)$ is bounded by a constant depending on $n$ only.
\end{enumerate}

Take a vertex $x_{0}$ of $G$, and for each integer $i\geq 0$, let $X_{i}=\{x\in V(G):{\rm dist}_{G}(x_{0},x)=i\}$.
Let $d=\max\{i:i\geq 0,~X_{i}\neq \emptyset \}~(={\rm ecc}_{G}(x_{0}))$.
Then $X_{0}=\{x_{0}\}$, $X_{1}=N_{G}(x_{0})$, and $V(G)$ is the disjoint union of $X_{i}~(0\leq i\leq d)$.
In this proof, we frequently use the following facts without mentioning:
\begin{enumerate}[{$\bullet $}]
\setlength{\parskip}{0cm}
\setlength{\itemsep}{0cm}
\item
For an integer $i\geq 1$ and a vertex $x\in X_{i}$, $N_{G}(x)\subseteq X_{i-1}\cup X_{i}\cup X_{i+1}$.
\item
For an integer $i\geq 0$ and a vertex $x\in X_{i}$, every shortest $x_{0}$-$x$ path $P$ of $G$ satisfies $|V(P)\cap X_{j}|=1$ for all $j$ with $0\leq j\leq i$.
\item
For integers $i$ and $i'$ with $0\leq i\leq i'$ and vertices $x\in X_{i}$ and $x'\in X_{i'}$, if an $x$-$x'$ path $P$ of $G$ satisfies $|V(P)\cap X_{j}|=1$ for all $j$ with $i\leq j\leq i'$, then $P$ is a shortest $x$-$x'$ path of $G$.
\end{enumerate}

\begin{claim}
\label{cl2.1}
Let $i$ and $k$ be integers with $n+1\leq i\leq k-n-1$.
Let $x\in X_{k}$, and let $Q$ be a shortest $x_{0}$-$x$ path of $G$.
Let $y\in X_{i}$ be a vertex with $N_{G}(y)\cap V(Q)\neq \emptyset $.
\begin{enumerate}[{\upshape(1)}]
\item
The vertex $y$ is adjacent to all vertices in $V(Q)\cap (X_{i-1}\cup X_{i+1})$.
\item
Let $i'$ be an integer with $n+1\leq i'\leq k-n-1$.
If a vertex $y'\in X_{i'}$ satisfies $N_{G}(y')\cap V(Q)=\emptyset $, then $yy'\notin E(G)$.
\end{enumerate}
\end{claim}
\proof
Write $Q=x_{0}x_{1}x_{2}\cdots x_{k}$ where $x_{k}=x$.
Note that $x_{i-n-1},x_{i-n},\ldots ,x_{i+n+1}$ are defined as vertices belonging to $V(Q)$.
\begin{enumerate}[{\upshape(1)}]
\item
We may assume that $y\notin V(Q)$.
If $N_{G}(y)\cap V(Q)=\{x_{j}\}$ for an integer $j$ with $j\in \{i-1,i,i+1\}$, then
$$
\{x_{j-n},x_{j-n+1},\ldots ,x_{j},y,x_{j+1},x_{j+2},\ldots ,x_{j+n}\}
$$
induces a copy of $F^{(1)}_{n}$ in $G$; if $N_{G}(y)\cap V(Q)=\{x_{j},x_{j+1}\}$ for an integer $j$ with $j\in \{i-1,i\}$, then
$$
\{x_{j-n+1},x_{j-n+2},\ldots ,x_{j},y,x_{j+1},x_{j+2},\ldots ,x_{j+n}\}
$$
induces a copy of $F^{(2)}_{n}$ in $G$.
In either case, we obtain a contradiction.
Considering the fact that $\emptyset \neq N_{G}(y)\cap V(Q)\subseteq \{x_{i-1},x_{i},x_{i+1}\}$, this forces $\{x_{i-1},x_{i+1}\}\subseteq N_{G}(y)\cap V(Q)$, as desired.
\item
Let $y'\in X_{i'}$ be a vertex with $N_{G}(y')\cap V(Q)=\emptyset $.
Note that $y'\notin V(Q)$.
Suppose that $yy'\in E(G)$.
By (1), $yx_{i-1},yx_{i+1}\in E(G)$, and hence $Q':=x_{0}x_{1}\cdots x_{i-1}yx_{i+1}x_{i+2}\cdots x_{k}$ is a shortest $x_{0}$-$x$ path of $G$.
Then $N_{G}(y')\cap V(Q')\supseteq \{y\}\neq \emptyset $.
Hence, applying (1) with $(y,Q)=(y',Q')$, $y'$ is adjacent to all vertices in $V(Q')\cap (X_{i'-1}\cup X_{i'+1})$.
Since $|V(Q')\cap (X_{i'-1}\cup X_{i'+1})|=2$ and $V(Q')\setminus V(Q)=\{y\}$, this implies that $y'$ is adjacent to a vertex in $V(Q)$, which is a contradiction.
\qed
\end{enumerate}

\begin{claim}
\label{cl2-rev-1}
Let $i$ and $i'$ be integers with $0\leq i<i'\leq d$.
\begin{enumerate}[{\upshape(1)}]
\item
There exists an induced star cover $\SS_{i,i'}$ of $G[\bigcup _{i\leq j\leq i'}X_{j}]$ with $|\SS_{i,i'}|\leq |X_{i}|c_{1}(n,2(i'-i))$, where $c_{1}(n,l_{0})$ is the constant appearing in Proposition~\ref{prop-indstar}.
\item
If $G$ is $\tilde{S}_{n}$-free, then there exists an induced star partition $\SS'_{i,i'}$ of $G[\bigcup _{i\leq j\leq i'}X_{j}]$ with $|\SS'_{i,i'}|\leq |X_{i}|c_{2}(n,2(i'-i))$, where $c_{2}(n,l_{0})$ is the constant appearing in Proposition~\ref{prop-indstar}.
\end{enumerate}
\end{claim}
\proof
Let $G'=G[\bigcup _{i\leq j\leq i'}X_{j}]$.
For an integer $j$ with $i+1\leq j\leq i'$ and a vertex $x\in X_{j}$, take a vertex $u_{x}\in X_{j-1}$ so that $u_{x}x\in E(G)$.
Let $D$ be the graph with $V(D)=\bigcup _{i\leq j\leq i'}X_{j}$ and $E(D)=\{u_{x}x:x\in X_{j},~i+1\leq j\leq i'\}$, and let $D_{1},D_{2},\ldots ,D_{q}$ be the components of $D$.

Fix an integer $h$ with $1\leq h\leq q$.
By the definition of $D$, $D_{h}$ is a tree and $|V(D_{h})\cap X_{i}|=1$, say $V(D_{h})\cap X_{i}=\{u_{h}\}$.
In particular, note that $q=|X_{i}|$.
Let $G_{h}=G[V(D_{h})]$.
Fix a vertex $x\in V(G_{h})$, and let $j$ be the integer with $i\leq j\leq i'$ and $x\in X_{j}$.
Then by the construction of $D$, we have ${\rm dist}_{G_{h}}(u_{h},x) = j-i\leq i'-i$.
Since $x$ is arbitrary, this implies that ${\rm diam}(G_{h})\leq 2{\rm ecc}_{G_{h}}(u_{h})\leq 2(i'-i)$.
Hence by applying Proposition~\ref{prop-indstar} with $l_{0}=2(i'-i)$,
\begin{align}
{\rm insc}(G_{h})\leq c_{1}(n,2(i'-i))\label{cond-cl2-rev-1-1}
\end{align}
and
\begin{align}
\mbox{if $G$ is $\tilde{S}_{n}$-free, then }{\rm insp}(G_{h})\leq c_{2}(n,2(i'-i)).\label{cond-cl2-rev-1-2}
\end{align}

Since $\bigcup _{i\leq j\leq i'}X_{i}~(=V(G'))$ is the disjoint union of $V(G_{h})~(1\leq h\leq q)$, it follows from (\ref{cond-cl2-rev-1-1}) and (\ref{cond-cl2-rev-1-2}) that
\begin{enumerate}[{$\bullet $}]
\setlength{\parskip}{0cm}
\setlength{\itemsep}{0cm}
\item
${\rm insc}(G')\leq \sum _{1\leq h\leq q}{\rm insc}(G_{h})\leq |X_{i}|c_{1}(n,2(i'-i))$, and
\item
if $G$ is $\tilde{S}_{n}$-free, then ${\rm insp}(G')\leq \sum _{1\leq h\leq q}{\rm insc}(G_{h})\leq |X_{i}|c_{2}(n,2(i'-i))$,
\end{enumerate}
as desired.
\qed

For the moment, we assume that $d\leq 4n+2$.
Note that the constants $c_{1}(n,8n+4)$ and $c_{2}(n,8n+4)$ are depending on $n$ only, where $c_{1}(n,l_{0})$ and $c_{2}(n,l_{0})$ are the constants appearing in Proposition~\ref{prop-indstar}.
By applying Claim~\ref{cl2-rev-1} with $i=0$ and $i'=d$,
$$
{\rm inspc}(G)\leq {\rm insc}(G)\leq c_{1}(n,2d)\leq c_{1}(n,8n+4),
$$
and
\begin{align*}
\mbox{if $G$ is $\tilde{S}_{n}$-free, then }{\rm inspp}(G)\leq {\rm insp}(G)\leq c_{2}(n,2d)\leq c_{2}(n,8n+4),
\end{align*}
and hence (B1) and (B2) hold.
Thus we may assume that
\begin{align}
d\geq 4n+3.\label{cond-sp-diam}
\end{align}

We consider the following operation recursively defining indices $k_{h}$, vertices $w_{h}$ and paths $Q_{h}$:
Let $k_{1}=d$, and take a vertex $w_{1}\in X_{k_{1}}$.
Let $Q_{1}$ be a shortest $x_{0}$-$w_{1}$ path of $G$.
Let $h$ be an integer with $h\geq 2$, and assume that $k_{1},\ldots ,k_{h-1}$, $w_{1},\ldots ,w_{h-1}$, $Q_{1},\ldots ,Q_{h-1}$ have been defined.
If $k_{h-1}\geq 3n+2$, let
$$
A_{h,i}=\left\{y\in X_{i}\setminus \left(\bigcup _{1\leq l\leq h-1}V(Q_{l})\right):N_{G}(y)\cap \left(\bigcup _{1\leq l\leq h-1}V(Q_{l})\right)=\emptyset \right\}
$$
for an integer $i$ with $2n+1\leq i\leq k_{h-1}-n-1$.
If $k_{h-1}\geq 3n+2$ and $A_{h,i}\neq \emptyset $ for some integer $i$ with $2n+1\leq i\leq k_{h-1}-n-1$, we let $k_{h}=\max\{i:2n+1\leq i\leq k_{h-1}-n-1,~A_{h,i}\neq \emptyset \}$ and $w_{h}\in A_{h,k_{h}}$, and let $Q_{h}$ be a shortest $x_{0}$-$w_{h}$ path of $G$; otherwise (i.e., either $k_{h-1}\geq 3n+2$ and $A_{h,i}=\emptyset $ for all integers $i$ with $2n+1\leq i\leq k_{h-1}-n-1$, or $k_{h-1}\leq 3n+1$), we let $k_{h}=2n$ and finish the operation.

Let $K=\{k_{h}:h\geq 1\}$, and set $h_{0}=|K|-1$.
Note that $d=k_{1}>k_{2}>\cdots >k_{h_{0}}>k_{h_{0}+1}=2n$.
For each integer $h$ with $1\leq h\leq h_{0}$, write $Q_{h}=v^{(h)}_{0}v^{(h)}_{1}v^{(h)}_{2}\cdots v^{(h)}_{k_{h}}$ where $v^{(h)}_{0}=x_{0}$ and $v^{(h)}_{k_{h}}=w_{h}$.
Note that $V(Q_{h})\cap X_{i}=\{v^{(h)}_{i}\}$ for integers $h$ and $i$ with $1\leq h\leq h_{0}$ and $1\leq i\leq k_{h}$.

\begin{figure}
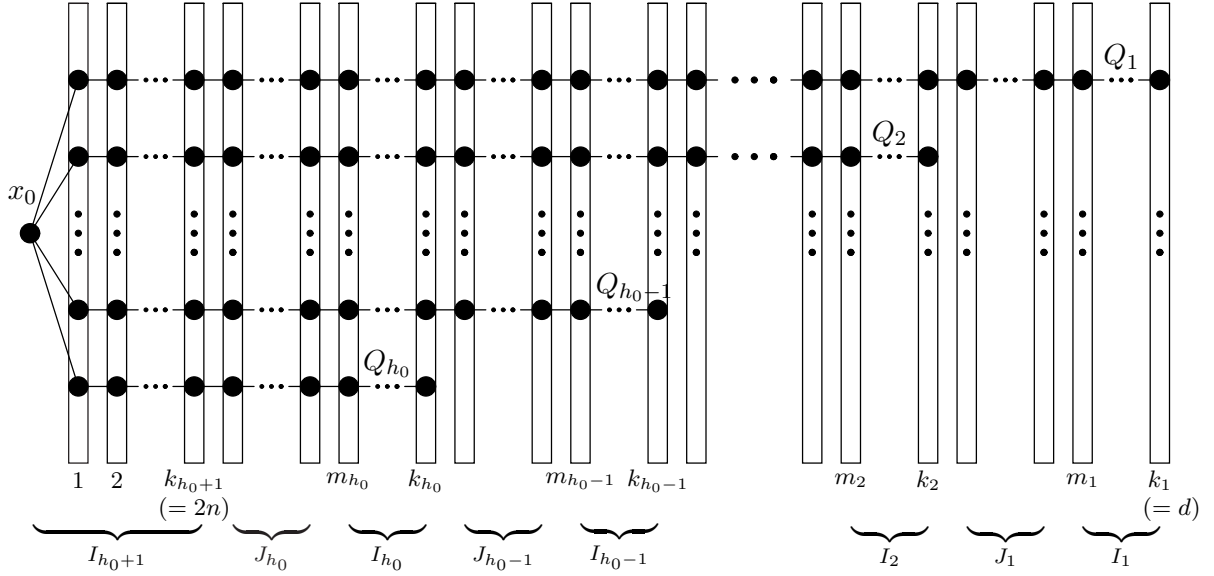

\hspace{-20mm}
{\unitlength 0.1in%
}%

\caption{Sets $I_{h}$ and $J_{h}$ of indices}
\label{f2}
\end{figure}

For an integer $h$ with $1\leq h\leq h_{0}$, we define the sets $I_{h}$ and $J_{h}$ of indices as
$$
I_{h}=
\begin{cases}
\{i:k_{h}-n\leq i\leq k_{h}\} & (1\leq h\leq h_{0}-1)\\
\{i:\max\{2n+1,k_{h_{0}}-n\}\leq i\leq k_{h_{0}}\} & (h=h_{0})
\end{cases}
$$
and
$$
J_{h}=\{i:k_{h+1}+1\leq i\leq k_{h}-n-1\}.
$$
Let $I_{h_{0}+1}=\{i:0\leq i\leq k_{h_{0}+1}~(=2n)\}$ (see Figure~\ref{f2}).
Note that $I_{h}$ is non-empty, but $J_{h}$ might be empty.
Since $k_{h_{0}+1}=2n$, we have $J_{h_{0}}=\{i:2n+1\leq i\leq k_{h_{0}}-n-1\}$.
Hence
\begin{align}
\mbox{$J_{h_{0}}=\emptyset $ if and only if $2n+1\geq k_{h_{0}}-n$.}\label{cl2.3--iff-1}
\end{align}

\begin{claim}
\label{cl2.3-}
Let $h$ be an integer with $1\leq h\leq h_{0}$, and let $i\in J_{h}$.
Then for every vertex $y\in X_{i}$, $N_{G}(y)\cap (\bigcup _{1\leq l\leq h}V(Q_{l}))\neq \emptyset $.
\end{claim}
\proof
By the definition of $k_{h}$, the desired conclusion clearly holds.
\qed

\begin{claim}
\label{cl2.3--}
The set $\{0,1,\ldots ,d\}$ is the disjoint union of $I_{h}~(1\leq h\leq h_{0}+1)$ and $J_{h}~(1\leq h\leq h_{0})$.
\end{claim}
\proof
We first prove that
\begin{align}
\mbox{$\{i:k_{h_{0}+1}+1\leq i\leq k_{h_{0}}\}$ is the disjoint union of $I_{h_{0}}$ and $J_{h_{0}}$.}\label{cond-cl2.3---1}
\end{align}
If $J_{h_{0}}=\emptyset $, then by (\ref{cl2.3--iff-1}), $\max\{2n+1,k_{h_{0}}-n\}=2n+1=k_{h_{0}+1}+1$, and hence $I_{h_{0}}=\{i:k_{h_{0}+1}+1\leq i\leq k_{h_{0}}\}$, as desired.
Thus we may assume that $J_{h_{0}}\neq \emptyset $.
Then by (\ref{cl2.3--iff-1}), $2n+1\leq k_{h_{0}}-n-1$, and hence $I_{h_{0}}=\{i:k_{h_{0}}-n\leq i\leq k_{h_{0}}\}$.
Since $J_{h_{0}}=\{i:k_{h_{0}+1}+1\leq i\leq k_{h_{0}}-n-1\}$, (\ref{cond-cl2.3---1}) holds.

By (\ref{cond-cl2.3---1}) and the definition of $I_{h}$ and $J_{h}$, the desired conclusion holds.
\qed

For an integer $h$ with $1\leq h\leq h_{0}$, let $m_{h}=\min I_{h}$.
Then $I_{h}=\{i:m_{h}\leq i\leq k_{h}\}$ and $J_{h}=\{i:k_{h+1}+1\leq i\leq m_{h}-1\}$ for every integer $h$ with $1\leq h\leq h_{0}$.




\begin{claim}
\label{cl2.3}
\begin{enumerate}[{\upshape(1)}]
\item
For two integers $h$ and $h'$ with $1\leq h<h'\leq h_{0}$, $V(Q_{h})\cap V(Q_{h'})=\{x_{0}\}$ and there is no edge between $V(Q_{h})\setminus \{x_{0}\}$ and $V(Q_{h'})\setminus \{x_{0}\}$.
\item
We have $h_{0}\leq 2$.
\item
For an integer $h$ with $1\leq h\leq h_{0}$, $|I_{h}|\leq n+1$.
\item
We have $\sum _{1\leq h\leq h_{0}+1}|I_{h}|\leq 4n+3$.
\end{enumerate}
\end{claim}
\proof
\begin{enumerate}[{\upshape(1)}]
\item
Since $x_{0}\in N_{G}(v^{(h')}_{1})\cap V(Q_{h})$, the value $p=\max\{i:1\leq i\leq k_{h'},~N_{G}(v^{(h')}_{i})\cap V(Q_{h})\neq \emptyset \}$ is well-defined.
By the definition of $w_{h'}$, we have $w_{h'}\in X_{k_{h'}}\setminus V(Q_{h})$ and $N_{G}(w_{h'})\cap V(Q_{h})=\emptyset $.
In particular, we have $1\leq p\leq k_{h'}-1$.
If there exists an integer $i$ with $p\leq i\leq k_{h'}-1$ such that $v^{(h')}_{i}\in V(Q_{h})$, then $N_{G}(v^{(h')}_{i+1})\cap V(Q_{h})\supseteq \{v^{(h')}_{i}\}\neq \emptyset $, which contradict the maximality of $p$.
Thus
\begin{align}
V(Q_{h})\cap \{v^{(h')}_{i}:p\leq i\leq k_{h'}\}=\emptyset .\label{cond-cl2.3-rev}
\end{align}

Suppose that $p\geq n+1$.
Recall that $N_{G}(v^{(h')}_{p})\cap V(Q_{h})\neq \emptyset $ and $N_{G}(v^{(h')}_{p+1})\cap V(Q_{h})=\emptyset $.
Since $n+1\leq p<p+1\leq k_{h'}\leq k_{h'-1}-n-1\leq k_{h}-n-1$, it follows from Claim~\ref{cl2.1}(2) with $(Q,y,y')=(Q_{h},v^{(h')}_{p},v^{(h')}_{p+1})$ that $v^{(h')}_{p}v^{(h')}_{p+1}\notin E(G)$, which is a contradiction.
Thus $p\leq n$.

Let $q=\max\{i:v^{(h')}_{p}v^{(h)}_{i}\in E(G)\}$.
Then $q\in \{p-1,p,p+1\}$.
Suppose that $q\geq 1$.
Since $p\leq n$ and $2n+1\leq k_{h'}$, we have
$$
q+n\leq p+n+1\leq 2n+1\leq k_{h'}<k_{h}.
$$
Hence $v^{(h)}_{q-1},v^{(h)}_{q},\ldots ,v^{(h)}_{q+n}$ have been defined as $n+2$ vertices belonging to $V(Q_{h})$, and $v^{(h')}_{p},v^{(h')}_{p+1},\ldots ,v^{(h')}_{p+n-1}$ have been defined as $n$ vertices belonging to $V(Q_{h'})$.
By (\ref{cond-cl2.3-rev}), $v^{(h')}_{p}\neq v^{(h)}_{q-1}$.
If $v^{(h')}_{p}v^{(h)}_{q-1}\notin E(G)$, then
$$
\{v^{(h')}_{p+n-1},v^{(h')}_{p+n-2},\ldots ,v^{(h')}_{p},v^{(h)}_{q},v^{(h)}_{q-1},v^{(h)}_{q+1},v^{(h)}_{q+2},\ldots ,v^{(h)}_{q+n}\}
$$
induces a copy of $F^{(1)}_{n}$ in $G$; if $v^{(h')}_{p}v^{(h)}_{q-1}\in E(G)$, then
$$
\{v^{(h')}_{p+n-1},v^{(h')}_{p+n-2},\ldots ,v^{(h')}_{p},v^{(h)}_{q-1},v^{(h)}_{q},\ldots ,v^{(h)}_{q+n-1}\}
$$
induces a copy of $F^{(2)}_{n}$ in $G$.
In either case, we obtain a contradiction.
Consequently, $q=0$, and this forces $p=1$.
Therefore $V(Q_{h})\cap V(Q_{h'})=\{x_{0}\}$ and there is no edge between $V(Q_{h})\setminus \{x_{0}\}$ and $V(Q_{h'})\setminus \{x_{0}\}$.

\item
It follows from (1) that if $h_{0}\geq 3$, then
$$
\{v^{(1)}_{n},v^{(1)}_{n-1},\ldots ,v^{(1)}_{1},x_{0},v^{(2)}_{1},v^{(3)}_{1},v^{(3)}_{2},\ldots ,v^{(3)}_{n}\}
$$
induces a copy of $F^{(1)}_{n}$ in $G$, which is a contradiction.

\item
By the definition of $I_{h}$ and (2),
\begin{enumerate}[{$\bullet $}]
\setlength{\parskip}{0cm}
\setlength{\itemsep}{0cm}
\item
$|I_{h_{0}}|=k_{h_{0}}-\max\{2n+1,k_{h_{0}}-n\}+1\leq k_{h_{0}}-(k_{h_{0}}-n)+1=n+1$, and
\item
$|I_{1}|=k_{1}-(k_{1}-n)+1=n+1$ if $h_{0}=2$,
\end{enumerate}
which leads to the desired conclusion.

\item
Since $|I_{h_{0}+1}|=2n+1$, it follows from (2) and (3) that
$$
\sum _{1\leq h\leq h_{0}+1}|I_{h}|\leq h_{0}(n+1)+2n+1\leq 2(n+1)+2n+1=4n+3,
$$
as desired.
\qed
\end{enumerate}

For integers $l$ with $1\leq l\leq h_{0}$ and $i\in \bigcup _{l\leq h\leq h_{0}}J_{h}$, let
$$
X^{(l)}_{i}=\{y\in X_{i}:N_{G}(y)\cap V(Q_{l})\neq \emptyset \}.
$$
Note that $v^{(l)}_{i}\in X^{(l)}_{i}$, and in particular, $X^{(l)}_{i}\neq \emptyset $.
Let $\nu =R(n-1,n)-1$.
Then $\nu $ is a constant depending on $n$ only.

\begin{claim}
\label{cl2.4}
Let $h$ and $i$ be integers with $1\leq h\leq h_{0}$ and $i\in J_{h}$.
\begin{enumerate}[{\upshape(1)}]
\item
The set $X_{i}$ is the disjoint union of $X^{(1)}_{i},X^{(2)}_{i},\ldots ,X^{(h)}_{i}$.
\item
If $i+1\in J_{h}$, then for each integer $l$ with $1\leq l\leq h$, every vertex in $X^{(l)}_{i}$ is adjacent to all vertices in $X^{(l)}_{i+1}$.
\item
For an integer $l$ with $1\leq l\leq h$, $|X^{(l)}_{i}|\leq \nu $.
\item
We have $|X_{i}|\leq h\nu $.
\item
If $G$ is $\{F^{(4)}_{n},F^{(5)}_{n}\}$-free, then $X_{i}\subseteq \bigcup _{1\leq l\leq h}V(Q_{l})$.
\end{enumerate}
\end{claim}
\proof
For an integer $l$ with $1\leq l\leq h$,
\begin{align}
n+1<2n+1\leq i\leq k_{h}-n-1\leq k_{l}-n-1,\label{eq-cl2.4-condition1}
\end{align}
and
\begin{align}
\mbox{if }i+1\in J_{h},\mbox{ then }i+1\leq k_{h}-n-1\leq k_{l}-n-1.\label{eq-cl2.4-condition2}
\end{align}

\begin{enumerate}[{\upshape(1)}]
\item
For a vertex $y\in X_{i}$, it follows from Claim~\ref{cl2.3-} that there exists an integer $l$ with $1\leq l\leq h$ such that $y\in X^{(l)}_{i}$.
This implies that $X_{i}=\bigcup _{1\leq l\leq h}X^{(l)}_{i}$.
Thus, to complete the proof of (1), it suffices to show that $X^{(1)}_{i},X^{(2)}_{i},\ldots ,X^{(h)}_{i}$ are pairwise disjoint.

Let $l_{1}$ and $l_{2}$ be integers with $1\leq l_{1}<l_{2}\leq h$.
Suppose that there exists a vertex $y\in X^{(l_{1})}_{i}\cap X^{(l_{2})}_{i}$.
Then $N_{G}(y)\cap V(Q_{l_{1}})\neq \emptyset $ and there exists a vertex $y'\in V(Q_{l_{2}})$ with $yy'\in E(G)$.
By Claim~\ref{cl2.3}(1), $N_{G}(y')\cap V(Q_{l_{1}})=\emptyset $.
Let $i'\in \{i-1,i,i+1\}$ be the integer with $y'\in X_{i'}$.
By (\ref{eq-cl2.4-condition1}), $n+1<i\leq k_{l_{2}}-n-1<k_{l_{1}}-n-1$, and hence $n+1\leq i-1\leq i'\leq i+1\leq k_{l_{1}}-n-1$.
Then by Claim~\ref{cl2.1}(2) with $Q=Q_{l_{1}}$, $yy'\notin E(G)$, which is a contradiction.
Thus $X^{(l_{1})}_{i}\cap X^{(l_{2})}_{i}=\emptyset $.
Since $l_{1}$ and $l_{2}$ are arbitrary, $X^{(1)}_{i},X^{(2)}_{i},\ldots ,X^{(h)}_{i}$ are pairwise disjoint.

\item
Fix two vertices $y\in X^{(l)}_{i}$ and $y'\in X^{(l)}_{i+1}$.
Since $N_{G}(y)\cap V(Q_{l})\neq \emptyset $ and $N_{G}(y')\cap V(Q_{l})\neq \emptyset $, it follows from (\ref{eq-cl2.4-condition1}), (\ref{eq-cl2.4-condition2}) and Claim~\ref{cl2.1}(1) with $Q=Q_{l}$ and $(y,Q)=(y',Q_{l})$ that $yv^{(l)}_{i-1},yv^{(l)}_{i+1},y'v^{(l)}_{i+2}\in E(G)$.
Hence $Q':=v^{(l)}_{0}v^{(l)}_{1}\cdots v^{(l)}_{i-1}yv^{(l)}_{i+1}v^{(l)}_{i+2}\cdots v^{(l)}_{k_{l}}$ is a shortest $x_{0}$-$w_{l}$ path of $G$ and $N_{G}(y')\cap V(Q')\neq \emptyset $.
Thus, again by (\ref{eq-cl2.4-condition2}) and Claim~\ref{cl2.1}(1) with $(y,Q)=(y',Q')$, $y'$ is adjacent to a vertex in $V(Q')\cap X_{i}~(=\{y\})$.
Since $y$ and $y'$ are arbitrary, the desired conclusion holds.

\item
It follows from (\ref{eq-cl2.4-condition1}) and Claim~\ref{cl2.1}(1) with $Q=Q_{l}$ that every vertex in $X^{(l)}_{i}$ is adjacent to $v^{(l)}_{i-1}$ in $G$, and hence $|X^{(l)}_{i}|\leq |N_{G}(v^{(l)}_{i-1})\cap X_{i}|$.
Thus, to complete the proof of (3), it suffices to show that $|N_{G}(v^{(l)}_{i-1})\cap X_{i}|\leq \nu $.
By way of contradiction, suppose that $|N_{G}(v^{(l)}_{i-1})\cap X_{i}|\geq \nu +1~(=R(n-1,n))$.
If there exists a clique $U\subseteq N_{G}(v^{(l)}_{i-1})\cap X_{i}$ of $G$ with $|U|=n-1$, then $\{v^{(l)}_{i-1}\}\cup U$ induces a copy of $K_{n}$ in $G$, which is a contradiction.
This together with the definition of the Ramsey number implies that there exists an independent set $U'\subseteq N_{G}(v^{(l)}_{i-1})\cap X_{i}$ of $G$ with $|U'|=n$.
Then again by (\ref{eq-cl2.4-condition1}) and Claim~\ref{cl2.1}(1) with $Q=Q_{l}$, every vertex in $U'$ is adjacent to $v^{(l)}_{i+1}$ in $G$.
Hence
$$
\{v^{(l)}_{i-n},v^{(l)}_{i-n+1},\ldots ,v^{(l)}_{i-1}\}\cup U'\cup \{v^{(l)}_{i+1},v^{(l)}_{i+2},\ldots ,v^{(l)}_{i+n}\}
$$
induces a copy of $F^{(3)}_{n}$ in $G$, which is a contradiction.

\item
By (1) and (3), $|X_{i}|=\sum _{1\leq l\leq h}|X^{(l)}_{i}|\leq h\nu $.


\item
Suppose that there exists a vertex $y\in X_{i}\setminus (\bigcup _{1\leq l\leq h}V(Q_{l}))$.
By (1), there exists an integer $l_{1}$ with $1\leq l_{1}\leq h$ such that $y\in X^{(l_{1})}_{i}$.
Then $N_{G}(y)\cap V(Q_{l_{1}})\neq \emptyset $.
It follows from (\ref{eq-cl2.4-condition1}) and Claim~\ref{cl2.1}(1) that $yv^{(l_{1})}_{i-1},yv^{(l_{1})}_{i+1}\in E(G)$, and hence
$$
\{v^{(l_{1})}_{i-n},v^{(l_{1})}_{i-n+1},\ldots ,v^{(l_{1})}_{i},y,v^{(l_{1})}_{i+1},v^{(l_{1})}_{i+2},\ldots ,v^{(l_{1})}_{i+n}\}
$$
induces a copy of $F^{(4)}_{n}$ or $F^{(5)}_{n}$ in $G$ according as $yv^{(l_{1})}_{i}\notin E(G)$ or not, which is a contradiction.
Thus $X_{i}\subseteq \bigcup _{1\leq l\leq h}V(Q_{l})$.
\qed
\end{enumerate}

For an integer $h$ with $1\leq h\leq h_{0}$, let $J'_{h}=J_{h}\setminus \{m_{h}-1\}~(=\{i:k_{h+1}+1\leq i\leq k_{h}-n-2\})$.

\begin{claim}
\label{cl2.4+}
Let $h$ be an integer with $1\leq h\leq h_{0}$.
\begin{enumerate}[{\upshape(1)}]
\item
There exists an induced path cover $\RR_{h}$ of $G[\bigcup _{i\in J'_{h}}X_{i}]$ with $|\RR_{h}|\leq h\nu $.
\item
If $G$ is $\{F^{(4)}_{n},F^{(5)}_{n}\}$-free, then there exists an induced path partition $\RR'_{h}$ of $G[\bigcup _{i\in J'_{h}}X_{i}]$ with $|\RR'_{h}|\leq h$.
\end{enumerate}
\end{claim}
\proof
If $J'_{h}=\emptyset $, then the claim clearly holds.
Thus we may assume that $J'_{h}\neq \emptyset $.
\begin{enumerate}[{\upshape(1)}]
\item
Fix an integer $l$ with $1\leq l\leq h$.
For $i\in J'_{h}$, write $X^{(l)}_{i}=\{a^{(l)}_{i,1},a^{(l)}_{i,2},\ldots ,a^{(l)}_{i,p_{l,i}}\}$ where $p_{l,i}=|X^{(l)}_{i}|$.
By Claim~\ref{cl2.4}(2), for $i\in J'_{h}\setminus \{m_{h}-2\}$, every vertex in $X^{(l)}_{i}$ is adjacent to all vertices in $X^{(l)}_{i+1}$.
Hence for an integer $j$ with $1\leq j\leq \nu $, $R^{(l)}_{j}:=b^{(l)}_{k_{h+1}+1,j}b^{(l)}_{k_{h+1}+2,j}\cdots b^{(l)}_{k_{h}-n-2,j}$ is an induced path of $G$ where
$$
b^{(l)}_{i,j}=
\begin{cases}
a^{(l)}_{i,j} & (j\leq p_{l,i})\\
a^{(l)}_{i,p_{l,i}} & (j>p_{l,i}).
\end{cases}
$$
By Claim~\ref{cl2.4}(3), $\max\{p_{l,i}:i\in J'_{h}\}=\max\{|X^{(l)}_{i}|:i\in J'_{h}\}\leq \nu $, and hence $(\bigcup _{1\leq j\leq \nu }V(R^{(l)}_{j}))\cap X_{i}=X^{(l)}_{i}$.
Consequently, it follows from Claim~\ref{cl2.4}(1) that $\RR_{h}:=\{R^{(l)}_{j}:1\leq l\leq h,~1\leq j\leq \nu \}$ is an induced path cover of $G[\bigcup _{i\in J'_{h}}(\bigcup _{1\leq l\leq h}X^{(l)}_{i})]~(=G[\bigcup _{i\in J'_{h}}X_{i}])$ with $|\RR_{h}|\leq h\nu $.

\item
By Claim~\ref{cl2.4}(5), $\bigcup _{i\in J'_{h}}X_{i}\subseteq \bigcup _{1\leq l\leq h}V(Q_{l})$.
Hence $\RR'_{h}:=\{Q_{l}[\bigcup _{i\in J'_{h}}X_{i}]:1\leq l\leq h\}$ is an induced path partition of $G[\bigcup _{i\in J'_{h}}X_{i}]$ with $|\RR'_{h}|=h$.
\qed
\end{enumerate}

Let $L=\{h:1\leq h\leq h_{0},~J_{h}\neq \emptyset \}$.
If $L=\emptyset $, then by Claims~\ref{cl2.3--} and \ref{cl2.3}(4),
$$
d+1 = |\{0,1,\ldots ,d\}|=\sum _{1\leq h\leq h_{0}+1}|I_{h}|\leq 4n+3,
$$
which contradicts (\ref{cond-sp-diam}).
Thus $L\neq \emptyset $.

Write $L=\{p_{1},\ldots ,p_{r}\}$ where $|L|=r$ and $p_{1}<\cdots <p_{r}$, and let $p_{0}=0$.
For an integer $h$ with $1\leq h\leq r+1$, let
$$
I'_{h}=
\begin{cases}
(\bigcup _{p_{h-1}+1\leq l\leq p_{h}}I_{l})\cup \{m_{p_{h}}-1\}~(=\{i:m_{p_{h}}-1\leq i\leq k_{p_{h-1}+1}\}) & (1\leq h\leq r)\\
\bigcup _{p_{r}+1\leq l\leq h_{0}+1}I_{l}~(=\{i:0\leq i\leq k_{p_{r}+1}\}) & (h=r+1).
\end{cases}
$$
Recall that $J'_{h}=J_{h}\setminus \{m_{h}-1\}$ for each $h$ with $1\leq h\leq h_{0}$.
Hence by Claim~\ref{cl2.3--},
\begin{align}
\mbox{$\{0,1,\ldots ,d\}$ is the disjoint union of $I'_{h}~(1\leq h\leq r+1)$ and $J'_{p_{h}}~(1\leq h\leq r)$.}\label{cond-sp-I'J'}
\end{align}


\begin{claim}
\label{cl2.5}
Let $c_{1}(n,l_{0})$ and $c_{2}(n,l_{0})$ be the constants appearing in Proposition~\ref{prop-indstar}.
\begin{enumerate}[{\upshape(1)}]
\item
For each integer $h$ with $1\leq h\leq r$, there exists an induced star cover $\SS_{h}$ of $G[\bigcup _{i\in I'_{h}}X_{i}]$ with $|\SS_{h}|\leq 2\nu c_{1}(n,4n+4)$.
\item
There exists an induced star cover $\SS_{r+1}$ of $G[\bigcup _{i\in I'_{r+1}}X_{i}]$ with $|\SS_{r+1}|\leq c_{1}(n,8n+4)$.
\item
For each integer $h$ with $1\leq h\leq r$, if $G$ is $\tilde{S}_{n}$-free, then there exists an induced star partition $\SS'_{h}$ of $G[\bigcup _{i\in I'_{h}}X_{i}]$ with $|\SS'_{h}|\leq 2\nu c_{2}(n,4n+4)$.
\item
If $G$ is $\tilde{S}_{n}$-free, then there exists an induced star partition $\SS'_{r+1}$ of $G[\bigcup _{i\in I'_{r+1}}X_{i}]$ with $|\SS'_{r+1}|\leq c_{2}(n,8n+4)$.
\end{enumerate}
\end{claim}
\proof
For each integer $h$ with $1\leq h\leq r+1$, let $G'_{h}=G[\bigcup _{i\in I'_{h}}X_{i}]$.

We first prove (2) and (4).
Since $I'_{r+1}=\bigcup _{p_{r}+1\leq l\leq h_{0}+1}I_{l}$, it follows from Claim~\ref{cl2.3}(4) that $|I'_{r+1}|\leq \sum _{1\leq l\leq h_{0}+1}|I_{l}|\leq 4n+3$.
Since $|X_{0}|=1$, this together with Claim~\ref{cl2-rev-1} implies that
$$
{\rm insc}(G'_{r+1})\leq |X_{0}|c_{1}(n,2(|I'_{r+1}|-1))\leq c_{1}(n,8n+4)
$$
and
\begin{align*}
\mbox{if $G$ is $\tilde{S}_{n}$-free, then }{\rm insp}(G'_{r+1})\leq |X_{0}|c_{2}(n,2(|I'_{r+1}|-1))\leq c_{2}(n,8n+4).
\end{align*}
Consequently, (2) and (4) hold.

Next we prove (1) and (3).
Let $h$ be an integer with $1\leq h\leq r$.
Since $I'_{h}=(\bigcup _{p_{h-1}+1\leq l\leq p_{h}}I_{l})\cup \{m_{p_{h}}-1\}$, it follows from Claim~\ref{cl2.3}(2)(3) that $|I'_{h}|-1\leq \sum _{p_{h-1}+1\leq l\leq p_{h}}|I_{l}|\leq \sum _{1\leq l\leq h_{0}}|I_{l}|\leq 2(n+1)$.
Since $m_{p_{h}}-1\in J_{p_{h}}$, it follows from Claims~\ref{cl2.3}(2) and \ref{cl2.4}(4) that $|X_{m_{p_{h}}-1}|\leq p_{h}\nu \leq h_{0}\nu \leq 2\nu $.
Hence by Claim~\ref{cl2-rev-1},
$$
{\rm insc}(G'_{h})\leq |X_{m_{p_{h}}-1}|c_{1}(n,2(|I'_{h}|-1))\leq 2\nu c_{1}(n,4n+4)
$$
and
\begin{align*}
\mbox{if $G$ is $\tilde{S}_{n}$-free, then }{\rm insp}(G'_{h})\leq |X_{m_{p_{h}}-1}|c_{2}(n,2(|I'_{h}|-1))\leq 2\nu c_{2}(n,4n+4).
\end{align*}
Consequently, (1) and (3) hold.
\qed

Let $\RR_{h}$ and $\SS_{h}$ (and $\RR'_{h}$ and $\SS'_{h}$ if $G$ is $\{\tilde{S}_{n},F^{(4)}_{n},F^{(5)}_{n}\}$-free) be the families as in Claims~\ref{cl2.4+} and \ref{cl2.5}.
By (\ref{cond-sp-I'J'}),
\begin{enumerate}[{$\bullet $}]
\setlength{\parskip}{0cm}
\setlength{\itemsep}{0cm}
\item
the family $\RR:=(\bigcup _{1\leq h\leq r}\RR_{p_{h}})\cup (\bigcup _{1\leq h\leq r+1}\SS_{h})$ is an induced SP-cover of $G$, and
\item
if $G$ is $\{\tilde{S}_{n},F^{(4)}_{n},F^{(5)}_{n}\}$-free, then the family $\RR':=(\bigcup _{1\leq h\leq r}\RR'_{p_{h}})\cup (\bigcup _{1\leq h\leq r+1}\SS'_{h})$ is an induced SP-partition of $G$.
\end{enumerate}
By Claim~\ref{cl2.3}(2), $r\leq h_{0}\leq 2$ and $1\leq p_{1}<\cdots <p_{r}\leq h_{0}\leq 2$, and in particular, $\sum _{1\leq h\leq r}p_{h}\leq 3$.
Hence
\begin{align*}
{\rm inspc}(G) &\leq |\RR|\\
&= \sum _{1\leq h\leq r}|\RR_{p_{h}}|+\sum _{1\leq h\leq r}|\SS_{h}|+|\SS_{r+1}|\\
&\leq \sum _{1\leq h\leq r}p_{h}\nu +2r\nu c_{1}(n,4n+4)+c_{1}(n,8n+4)\\
&\leq 3\nu +4\nu c_{1}(n,4n+4)+c_{1}(n,8n+4),
\end{align*}
and if $G$ is $\{\tilde{S}_{n},F^{(4)}_{n},F^{(5)}_{n}\}$-free, then
\begin{align*}
{\rm inspp}(G) &\leq |\RR'|\\
&= \sum _{1\leq h\leq r}|\RR'_{p_{h}}|+\sum _{1\leq h\leq r}|\SS'_{h}|+|\SS'_{r+1}|\\
&\leq \sum _{1\leq h\leq r}p_{h}+2r\nu c_{2}(n,4n+4)+c_{2}(n,8n+4)\\
&\leq 3+4\nu c_{2}(n,4n+4)+c_{2}(n,8n+4).
\end{align*}
This completes the proof of Proposition~\ref{mainthm-2}.
\qed

We conclude this section by finding families $\HH$ satisfying {\rm (P-${\rm insc}$)}, {\rm (P-${\rm insp}$)}, {\rm (P-${\rm inpc}$)} and {\rm (P-${\rm inpp}$)}.

\begin{prop}
\label{mainthm-star/path}
Let $n\geq 4$ be an integer.
Then the following hold:
\begin{enumerate}
\item[{\upshape(i)}]
There exists a constant $c_{\rm insc}(n)$ depending on $n$ only such that ${\rm insc}(G)\leq c_{\rm insc}(n)$ for every connected $\{K_{n},S^{*}_{n},P_{n}\}$-free graph $G$.
\item[{\upshape(ii)}]
There exists a constant $c_{\rm insp}(n)$ depending on $n$ only such that ${\rm insp}(G)\leq c_{\rm insp}(n)$ for every connected $\{K_{n},S^{*}_{n},\tilde{S}_{n},P_{n}\}$-free graph $G$.
\item[{\upshape(iii)}]
There exists a constant $c_{\rm inpc}(n)$ depending on $n$ only such that ${\rm inpc}(G)\leq c_{\rm inpc}(n)$ for every connected $\{K_{n},K_{1,n},F^{(1)}_{n},F^{(2)}_{n}\}$-free graph $G$.
\item[{\upshape(iv)}]
There exists a constant $c_{\rm inpp}(n)$ depending on $n$ only such that ${\rm inpp}(G)\leq c_{\rm inpp}(n)$ for every connected $\{K_{n},K_{1,n},F^{(1)}_{n},F^{(2)}_{n},F^{(4)}_{n},F^{(5)}_{n}\}$-free graph $G$.
\end{enumerate}
\end{prop}
\proof
Let $c_{\rm inspc}(n)$ and $c_{\rm inspp}(n)$ be the constants appearing in the statement in Proposition~\ref{mainthm-2}.

Let $G$ be a connected $\{K_{n},S^{*}_{n},P_{n}\}$-free graph.
By the $P_{n}$-freeness of $G$, $G$ is also $\{F^{(1)}_{n},F^{(2)}_{n},F^{(3)}_{n}\}$-free.
Hence by Proposition~\ref{mainthm-2}(i), there exists an induced SP-cover $\RR$ of $G$ with $|\RR|\leq c_{\rm inspc}(n)$.
Let $\PP=\{P\in \RR:P\mbox{ is a path}\}$, and let $\SS=\{\{x\}:x\in V(P),~P\in \PP\}$.
Then $(\RR\setminus \PP)\cup \SS$ is an induced star cover of $G$.
Since $G$ is $P_{n}$-free, $|V(P)|\leq n-1$ for every $P\in \PP$.
Hence $|\SS|=|\bigcup _{P\in \PP}V(P)|\leq \sum _{P\in \PP}|V(P)|\leq (n-1)|\PP|\leq (n-1)c_{\rm inspc}(n)$.
Therefore,
$$
{\rm insc}(G)\leq |\RR\setminus \PP|+|\SS|\leq c_{\rm inspc}(n)+(n-1)c_{\rm inspc}(n).
$$
Since $G$ is arbitrary, (i) holds.

By similar argument, we see that ${\rm insp}(G)\leq nc_{\rm inspp}(n)$ for every connected $\{K_{n},S^{*}_{n},\tilde{S}_{n},P_{n}\}$-free graph $G$.
Furthermore, since every $K_{1,n}$-free graph is also $\{S^{*}_{n},\tilde{S}_{n},F^{(3)}_{n}\}$-free and an induced star in a $K_{1,n}$-free graph has at most $n$ vertices, we also verify that
\begin{enumerate}[{$\bullet $}]
\setlength{\parskip}{0cm}
\setlength{\itemsep}{0cm}
\item
${\rm inpc}(G)\leq (n+1)c_{\rm inspc}(n)$ for every connected $\{K_{n},K_{1,n},F^{(1)}_{n},F^{(2)}_{n}\}$-free graph $G$, and
\item
${\rm inpp}(G)\leq (n+1)c_{\rm inspp}(n)$ for every connected $\{K_{n},K_{1,n},F^{(1)}_{n},F^{(2)}_{n},F^{(4)}_{n},F^{(5)}_{n}\}$-free graph $G$.
\end{enumerate}
Therefore, (ii)--(iv) also hold.
\qed

\section{Graphs with large induced cover/partition numbers}\label{sec-nec}

In this section, we calculate induced cover/partition numbers of some specific graphs.
The results will be used for proving the ``only if'' parts in main results.

We can easily verify that the following lemma holds, and so we omit its details.

\begin{lem}
\label{lem-nec-complete/star}
Let $n\geq 2$ be an integer.
Then ${\rm inspc}(K_{n})=\lceil \frac{n}{2} \rceil$, ${\rm inspc}(S^{*}_{n})=\lceil \frac{n}{2} \rceil $, ${\rm inspp}(\tilde{S}_{n})=n+1$, ${\rm inpc}(K_{1,n})=\lceil \frac{n}{2} \rceil$, and ${\rm insc}(P_{n})=\lceil \frac{n}{3} \rceil$.
\end{lem}

Let $m\geq 2$ and $n\geq 3$ be integers, and let $Q_{i}=u^{(1)}_{i}u^{(2)}_{i}\cdots u^{(n)}_{i}~(1\leq i\leq m)$ be $m$ pairwise vertex-disjoint paths.
We define five graphs.
\begin{enumerate}[{$\bullet $}]
\setlength{\parskip}{0cm}
\setlength{\itemsep}{0cm}
\item
Let $H^{(1)}_{m,n}$ be the graph obtained from the union of the paths $Q_{1},\ldots ,Q_{m}$ by adding $2(m-1)$ vertices $v_{i},w_{i}~(1\leq i\leq m-1)$ and $3(m-1)$ edges $v_{i}w_{i},v_{i}u^{(n)}_{i},v_{i}u^{(1)}_{i+1}~(1\leq i\leq m-1)$.
\item
Let $H^{(2)}_{m,n}$ be the graph obtained from the union of the paths $Q_{1},\ldots ,Q_{m}$ by adding $m-1$ vertices $v_{i}~(1\leq i\leq m-1)$ and $3(m-1)$ edges $v_{i}u^{(n)}_{i},v_{i}u^{(1)}_{i+1},u^{(n)}_{i}u^{(1)}_{i+1}~(1\leq i\leq m-1)$.
\item
Let $H^{(3)}_{m,n}$ be the graph obtained from the union of the paths $Q_{1},\ldots ,Q_{m}$ by adding $(m-1)m$ vertices $v_{i,j}~(1\leq i\leq m-1,~1\leq j\leq m)$ and $2(m-1)m$ edges $v_{i,j}u^{(n)}_{i},v_{i,j}u^{(1)}_{i+1}~(1\leq i\leq m-1,~1\leq j\leq m)$.
\item
Let $H^{(4)}_{m,n}$ be the graph obtained from $H^{(3)}_{m,n}$ by deleting $(m-1)(m-2)$ vertices $v_{i,j}~(1\leq i\leq m-1,~3\leq j\leq m)$.
\item
Let $H^{(5)}_{m,n}$ be the graph obtained from $H^{(4)}_{m,n}$ by adding $m-1$ edges $v_{i,1}v_{i,2}~(1\leq i\leq m-1)$.
\end{enumerate}
Then the following lemma holds.
(Here we just want to claim that ${\rm inspc}(H^{(i)}_{m,n})\rightarrow \infty $ when $m\rightarrow \infty $ for $i\in \{1,2,3\}$, and ${\rm inspp}(H^{(i)}_{m,n})\rightarrow \infty $ when $m\rightarrow \infty $ for $i\in \{4,5\}$.
Thus the readers recognizing the facts can skip the proof.)

\begin{lem}
\label{lem-nec-pathtype}
Let $m$ and $n$ be integers with $m\geq 2$ and $n\geq 3$.
\begin{enumerate}[{\upshape(1)}]
\item
We have ${\rm inspc}(H^{(1)}_{m,n})\geq \lceil \frac{m+1}{2} \rceil $.
\item
We have ${\rm inspc}(H^{(2)}_{m,n})\geq \lceil \frac{m+1}{2} \rceil $.
\item
We have ${\rm inspc}(H^{(3)}_{m,n})\geq m$.
\item
We have ${\rm inspp}(H^{(4)}_{m,n})\geq m$.
\item
We have ${\rm inspp}(H^{(5)}_{m,n})\geq m$.
\end{enumerate}
\end{lem}
\proof
\begin{enumerate}
\item[{\upshape(1)}]
Since any two leaves of $H^{(1)}_{m,n}$ have distance at least four, each element of an induced SP-cover of $H^{(1)}_{m,n}$ contains at most two leaves of $H^{(1)}_{m,n}$.
Since $H^{(1)}_{m,n}$ has $m+1$ leaves, this implies that ${\rm inspc}(H^{(1)}_{m,n})\geq \lceil \frac{m+1}{2} \rceil $.

\item[{\upshape(2)}]
Let $\PP$ be an induced SP-cover of $H^{(2)}_{m,n}$.
It suffices to show that $|\PP|\geq \lceil \frac{m+1}{2} \rceil $.
Since $H^{(2)}_{m,n}$ is $K_{1,3}$-free, $\PP$ can be regarded as an induced path cover of $H^{(2)}_{m,n}$.
Let $U=\{u^{(1)}_{1},u^{(n)}_{m}\}\cup \{v_{i}:1\leq i\leq m-1\}$ and ${\bf X}=\{(P,u)\in \PP\times U:u\in V(P)\}$.
For each $u\in U$, there exists an element $P$ of $\PP$ with $u\in V(P)$.
Hence
\begin{align}
|{\bf X}|=\sum _{u\in U}|\{P\in \PP:u\in V(P)\}|\geq |U|=m+1.\label{eq-lem-H2-1}
\end{align}
For every vertex $u\in U$, since the neighborhood of $u$ is a clique of $H^{(2)}_{m,n}$, no induced path of $H^{(2)}_{m,n}$ contains $u$ as an internal vertex.
Since a path contains at most two endvertices, this implies that $|V(P)\cap U|\leq 2$ for every $P\in \PP$.
Hence
\begin{align}
|{\bf X}|=\sum _{P\in \PP}|V(P)\cap U|\leq 2|\PP|.\label{eq-lem-H2-2}
\end{align}
By (\ref{eq-lem-H2-1}) and (\ref{eq-lem-H2-2}), we have $m+1\leq |{\bf X}|\leq 2|\PP|$, as desired.

\item[{\upshape(3)}]
Let $\PP$ be an induced SP-cover of $H^{(3)}_{m,n}$.
It suffices to show that $|\PP|\geq m$.
Let $U=\{u^{(1)}_{1}\}\cup \{v_{i,j}:1\leq i\leq m-1,~1\leq j\leq m\}$ and ${\bf X}=\{(P,u)\in \PP\times U:u\in V(P)\}$.
For each $u\in U$, there exists an element $P$ of $\PP$ with $u\in V(P)$.
Hence
\begin{align}
|{\bf X}|=\sum _{u\in U}|\{P\in \PP:u\in V(P)\}|\geq |U|=(m-1)m+1.\label{eq-lem-H3-1}
\end{align}
Note that for an integer $i$ with $1\leq i\leq m-1$, if an induced path $P$ of $H^{(3)}_{m,n}$ contains at least two vertices in $\{v_{i,j}:1\leq j\leq m\}$, then $|V(P)|=3$ and $P$ contains exactly two vertices in $\{v_{i,j}:1\leq j\leq m\}$.
Hence if an element $P$ of $\PP$ is a path, then $|V(P)\cap U|\leq m$.
(For example, a $u^{(1)}_{1}$-$u^{(n)}_{m}$ path $P$ of $H^{(3)}_{m,n}$ is an induced path of $H^{(3)}_{m,n}$ with $|V(P)\cap U|=m$.)
Since any two sets of $\{u^{(1)}_{1}\}$ and $\{v_{i,j}:1\leq j\leq m\}~(1\leq i\leq m-1)$ have distance at least three in $H^{(3)}_{m,n}$, if an element $P$ of $\PP$ is a star, then $|V(P)\cap U|\leq m$.
This implies that
\begin{align}
|{\bf X}|=\sum _{P\in \PP}|V(P)\cap U|\leq m|\PP|.\label{eq-lem-H3-2}
\end{align}
By (\ref{eq-lem-H3-1}) and (\ref{eq-lem-H3-2}), we have $m^{2}-m+1\leq |{\bf X}|\leq m|\PP|$, and hence $|\PP|\geq m-1+\frac{1}{m}$.
Since $|\PP|$ is an integer, this forces $|\PP|\geq m$.

\item[{\upshape(4)}]
Let $\PP$ be an induced SP-partition of $H^{(4)}_{m,n}$.
It suffices to show that $|\PP|\geq m$.
For an integer $i$ with $1\leq i\leq m$, let $R_{i}$ be the unique element of $\PP$ containing $u^{(2)}_{i}$.
We remark that $R_{i}$ might equal $R_{i'}$ for some $1\leq i<i'\leq m$.
Let $I=\{i:1\leq i\leq m-1,~R_{i}\neq R_{i+1}\}$, and write $I=\{i_{1},i_{2},\ldots ,i_{h}\}$ with $i_{1}<i_{2}<\ldots <i_{h}$ where $h=0$ if $I=\emptyset $.
If $R_{i_{l}}=R_{i'}$ for some integers $l$ and $i'$ with $1\leq l\leq h$ and $i'>i_{l}$, then $R_{i_{l}}$ is a path and contains all vertices in $\{u^{(2)}_{j}:i_{l}\leq j\leq i'\}$, and hence $R_{i_{l}}=R_{i_{l}+1}$, which contradicts the definition of $i_{l}$.
This, together with the fact that $i_{h}<m$ if $I\neq \emptyset $, implies that $R_{i_{1}},R_{i_{2}},\ldots ,R_{i_{h}},R_{m}$ are pairwise distinct, and so $|\{R_{i}:1\leq i\leq m\}|\geq h+1$.
Note that the inequality also holds even if $I=\emptyset $.

Fix an integer $j\in \{1,2,\ldots ,m-1\}\setminus I$.
By the definition of $I$, $R_{j}=R_{j+1}$, and hence $R_{j}$ contains $u^{(2)}_{j}$ and $u^{(2)}_{j+1}$.
Since ${\rm dist}_{H^{(4)}_{m,n}}(u^{(2)}_{j},u^{(2)}_{j+1})=n+1\geq 4$, $R_{j}$ is a path.
We can easily verify that $\{u^{(n)}_{j},u^{(1)}_{j+1}\}\subseteq V(R_{j})$ and $|V(R_{j})\cap \{v_{j,1},v_{j,2}\}|=1$.
This implies that there exists an element $R'_{j}$ of $\PP$ with $V(R'_{j})=\{v_{j,1}\}$ or $V(R'_{j})=\{v_{j,2}\}$.
Therefore
$$
|\PP|\geq |\{R_{i}:1\leq i\leq m\}\cup \{R'_{j}:j\in \{1,2,\ldots ,m-1\}\setminus I\}|\geq (h+1)+((m-1)-h)=m,
$$
as desired.

\item[{\upshape(5)}]
We proceed by induction on $m$.
Since $H^{(5)}_{2,n}$ contains a cycle, we have ${\rm inspp}(H^{(5)}_{2,n})\geq 2$.
Thus we may assume that $m\geq 3$.
We let $\PP $ be an induced SP-partition of $H^{(5)}_{m,n}$.
It suffices to show that $|\PP|\geq m$.
Since $H^{(5)}_{m,n}$ is $K_{1,3}$-free, $\PP$ can be regarded as an induced path partition of $H^{(5)}_{m,n}$.

For the moment, suppose that $\{v_{i,1},v_{i,2}\}\not\subseteq V(P)$ for any $P\in \PP$ and any integer $i$ with $1\leq i\leq m-1$.
Then $\PP$ is an induced SP-partition of $H^{(5)}_{m,n}-\{v_{i,1}v_{i,2}:1\leq i\leq m-1\}~(=H^{(4)}_{m,n})$.
Hence by (4), $|\PP|\geq {\rm inspp}(H^{(4)}_{m,n})\geq m$, as desired.
Thus we may assume that there exist an element $P^{*}$ of $\PP$ and an integer $i$ with $1\leq i\leq m-1$ such that $\{v_{i,1},v_{i,2}\}\subseteq V(P^{*})$.
Note that a vertex in $V(H^{(5)}_{m,n})\setminus \{v_{i,1},v_{i,2}\}$ is adjacent to $v_{i,1}$ if and only if the vertex is adjacent to $v_{i,2}$.
Since $P^{*}$ is an induced path of $H^{(5)}_{m,n}$, this implies that $P^{*}$ consists of two vertices $v_{i,1}$ and $v_{i,2}$.
In other words, $\PP\setminus \{P^{*}\}$ is an induced SP-partition of $H^{(5)}_{m,n}-\{v_{i,1},v_{i,2}\}$.
If $i\in \{1,m-1\}$, then $H^{(5)}_{m,n}-\{v_{i,1},v_{i,2}\}$ is the disjoint union of a path and a copy of $H^{(5)}_{m-1,n}$, and hence by the induction hypothesis, $|\PP\setminus \{P^{*}\}|\geq 1+{\rm inspp}(H^{(5)}_{m-1,n})\geq 1+(m-1)$; if $2\leq i\leq m-2$, then $H^{(5)}_{m,n}-\{v_{i,1},v_{i,2}\}$ is the disjoint union of a copy of $H^{(5)}_{i,n}$ and a copy of $H^{(5)}_{m-i,n}$, and hence $|\PP\setminus \{P^{*}\}|\geq {\rm inpp}(H^{(5)}_{i,n})+{\rm inpp}(H^{(5)}_{m-i,n})\geq i+(m-i)$.
In either case, we obtain $|\PP|>|\PP\setminus \{P^{*}\}|\geq m$, as desired.
\qed
\end{enumerate}

\section{Proof of main results}\label{sec-proof}

In this section, we complete the proof of Theorem~\ref{mainthm} and settle Ramsey-type problems for four invariants ${\rm insc}$, ${\rm insp}$, ${\rm inpc}$ and ${\rm inpp}$.

\medbreak\noindent\textit{Proof of Theorem~\ref{mainthm}.}\quad
By Proposition~\ref{mainthm-2}, we obtain the ``if'' parts of (i) and (ii).
Thus we prove that the ``only if'' parts are true.
Let $\HH$ be a finite family of connected graphs.
Then the value $p=\max\{|V(H)|:H\in \HH\}$ is a well-defined constant depending on $\HH$.

We first suppose that $\HH$ satisfies {\rm (P-${\rm inspc}$)}.
Then there exists a constant $c=c(\HH)$ such that ${\rm inspc}(G)\leq c$ for every connected $\HH$-free graph $G$.
Let $n=\max\{p,2c+1,4\}$.
Since ${\rm inspc}(K_{n})={\rm inspc}(S^{*}_{n})=\lceil \frac{n}{2} \rceil \geq c+1$ by Lemma~\ref{lem-nec-complete/star}, neither $K_{n}$ nor $S^{*}_{n}$ is $\HH$-free.
This implies that $\HH\leq \{K_{n},S^{*}_{n}\}$.
Furthermore, for each $i\in \{1,2,3\}$, it follows from Lemma~\ref{lem-nec-pathtype}(1)--(3) that ${\rm inspc}(H^{(i)}_{n,n})\geq c+1$, and hence $H^{(i)}_{n,n}$ is not $\HH$-free, i.e., $H^{(i)}_{n,n}$ contains an induced subgraph $A_{i}$ isomorphic to an element of $\HH$.
Since $A_{i}$ is connected and $|V(A_{i})|\leq \max\{|V(H)|:H\in \HH\}=p\leq n$, we have
\begin{enumerate}[$\bullet $]
\setlength{\parskip}{0cm}
\setlength{\itemsep}{0cm}
\item
$|\{i':1\leq i'\leq n,~V(A_{i})\cap V(Q_{i'})\neq \emptyset \}|\leq 2$,
\item
$|\{i':1\leq i'\leq n-1,~V(A_{1})\cap \{v_{i'},w_{i'}\}\neq \emptyset \}|\leq 1$,
\item
$|\{i':1\leq i'\leq n-1,~V(A_{2})\cap \{v_{i'}\}\neq \emptyset \}|\leq 1$, and
\item
$|\{i':1\leq i'\leq n-1,~V(A_{3})\cap \{v_{i',j}:1\leq j\leq n\}\neq \emptyset \}|\leq 1$.
\end{enumerate}
This implies that $A_{i}$ is isomorphic to an induced subgraph of $F^{(i)}_{n}$, and hence $\HH\leq \{F^{(1)}_{n},F^{(2)}_{n},F^{(3)}_{n}\}$.
Consequently, $\HH\leq \{K_{n},S^{*}_{n},F^{(1)}_{n},F^{(2)}_{n},F^{(3)}_{n}\}$, which proves the ``only if'' part of (i).

Next we suppose that $\HH$ satisfies {\rm (P-${\rm inspp}$)}.
Then there exists a constant $c'=c'(\HH)$ such that ${\rm inspp}(G)\leq c'$ for every connected $\HH$-free graph $G$.
Let $n=\max\{p,2c'+1,4\}$.
Note that for every graph $H$, ${\rm inspc}(H)\leq {\rm inspp}(H)$.
Hence by similar argument as above, together with Lemmas~\ref{lem-nec-complete/star} and \ref{lem-nec-pathtype}(1)(2)(4)(5), we see that $\HH\leq \{K_{n},S^{*}_{n},\tilde{S}_{n},F^{(1)}_{n},F^{(2)}_{n},F^{(4)}_{n},F^{(5)}_{n}\}$, which proves the ``only if'' part of (ii).

This completes the proof of Theorem~\ref{mainthm}.
\qed

Note that for every graph $H$,
\begin{align*}
&{\rm inspc}(H)\leq {\rm insc}(H)\leq {\rm insp}(H),~{\rm inspp}(H)\leq {\rm insp}(H),\\
&{\rm inspc}(H)\leq {\rm inpc}(H)\leq {\rm inpp}(H),\mbox{ and }{\rm inspp}(H)\leq {\rm inpp}(H).
\end{align*}
Hence, by similar argument in the proof of Theorem~\ref{mainthm}, together with Proposition~\ref{mainthm-star/path} and Lemmas~\ref{lem-nec-complete/star} and \ref{lem-nec-pathtype}(1)(2)(4)(5), we obtain the following theorem.

\begin{thm}
\label{mainthm-cor}
Let $\HH$ be a finite family of connected graphs.
\begin{enumerate}[{\upshape(i)}]
\item
The family $\HH$ satisfies {\rm (P-${\rm insc}$)} if and only if $\HH\leq \{K_{n},S^{*}_{n},P_{n}\}$ for an integer $n\geq 4$.
\item
The family $\HH$ satisfies {\rm (P-${\rm insp}$)} if and only if $\HH\leq \{K_{n},S^{*}_{n},\tilde{S}_{n},P_{n}\}$ for an integer $n\geq 4$.
\item
The family $\HH$ satisfies {\rm (P-${\rm inpc}$)} if and only if $\HH\leq \{K_{n},K_{1,n},F^{(1)}_{n},F^{(2)}_{n}\}$ for an integer $n\geq 4$.
\item
The family $\HH$ satisfies {\rm (P-${\rm inpp}$)} if and only if $\HH\leq \{K_{n},K_{1,n},F^{(1)}_{n},F^{(2)}_{n},F^{(4)}_{n},F^{(5)}_{n}\}$ for an integer $n\geq 4$.
\end{enumerate}
\end{thm}

Now we recall Remark~\ref{remark-mainthm-2} in Section~\ref{sec3}.
Let $G$ be a graph.
A family $\PP$ of subgraphs of $G$ is called an {\it isometric path cover} of $G$ if $\bigcup _{P\in \PP}V(P)=V(G)$ and each element of $\PP$ is an isometric path of $G$.
An isometric path cover $\PP$ of $G$ is called an {\it isometric path partition} of $G$ if the elements of $\PP$ are pairwise vertex-disjoint.
The minimum cardinality of an isometric path cover (resp. an isometric path partition) of $G$, denoted by ${\rm ispc}(G)$ (resp. ${\rm ispp}(G)$), is called the {\it isometric path cover number} (resp. the {\it isometric path partition number}) of $G$.
It is known that the isometric path cover number is closely related to the cops-robbers game, and motivated by the fact, there are many researches on the isometric path cover/partition (see \cite{FF,F2,FNHC,M,PC,PC2}).
The isometric path cover number is frequently called the {\it isometric path number}.
As we mentioned in Remark~\ref{remark-mainthm-2}, Proposition~\ref{mainthm-star/path}(iii)(iv) can be extended to isometric versions.
Furthermore, for every graph $G$, since every isometric path of $G$ is also an induced path of $G$, we have ${\rm inpc}(G)\leq {\rm ispc}(G)$ and ${\rm inpp}(G)\leq {\rm ispp}(G)$.
Consequently, we have implicitly proved the following theorem.

\begin{thm}
\label{mainthm-3}
Let $\HH$ be a finite family of connected graphs.
\begin{enumerate}[{\upshape(i)}]
\item
The family $\HH$ satisfies {\rm (P-${\rm ispc}$)} if and only if $\HH\leq \{K_{n},K_{1,n},F^{(1)}_{n},F^{(2)}_{n}\}$ for an integer $n\geq 4$.
\item
The family $\HH$ satisfies {\rm (P-${\rm ispp}$)} if and only if $\HH\leq \{K_{n},K_{1,n},F^{(1)}_{n},F^{(2)}_{n},F^{(4)}_{n},F^{(5)}_{n}\}$ for an integer $n\geq 4$.
\end{enumerate}
\end{thm}

\section*{Acknowledgment}

The authors would like to thank Professor Jin Sun for the comment improving the bound in Claim~\ref{cl2.3}(2) and an idea for shortening the proof of Lemma~\ref{lem-2-indstar}, and to thank the referees for their valuable suggestions and comments.
This work was partially supported by JSPS KAKENHI Grant number JP20K03720 (to S.C) and JSPS KAKENHI Grant number JP18K13449 and 23K03204 (to M.F).

\end{document}